\documentclass{jnmp}
\usepackage{amsmath,amssymb,amsthm}
\usepackage{array}
\usepackage{latexsym,epsfig}

\JNMPnumberwithin{equation}{section}

\newtheorem{coro}{Corollary}
\newtheorem{theo}{Proposition}

\newtheorem{rem}{Remark}

\newcommand{\sech}{\mathop{\mathrm{sech}}}
\newcommand{\Ref}[1]{{\normalfont(\ref{#1})}}
\newcommand{\AB}{\allowbreak}
\newcommand{\gr}[2]{\mathop{{\mathbf{#1}}(#2)}\nolimits}

\newcommand{\SU}[1]{\gr{SU}{#1}}
\newcommand{\ali}[2]{\mathop{\mathfrak{#1}(#2)}\nolimits}
\newcommand{\lie}[3]{\mathop{\mathfrak{#1}\,(#2, #3)}\nolimits}

\newcommand{\IS}[2]{\sum\limits_{#1}^{#2}}

\newcommand{\ID}{\mathop{\mathbf1}\nolimits}
\newcommand{\Tr}{\mathop{\mathrm{tr}}\nolimits}
\newcommand{\ad}{\mathop{\mathrm{ad}}\nolimits}
\newcommand{\ADA}[1]{\ifmmode \ad(#1) \else $\ad(#1)$\fi}
\newcommand{\LI}[2]{\ifmmode#2_1,\AB\,\ldots,\,\AB #2_{#1}%
\else$ #2_1,\AB\,\ldots,\,\AB#2_{#1}$\fi}

\newcommand{\su}[1]{\ali{su}{#1}}
\newcommand{\sltwo}{\ifmmode \ali{sl}{2} \else $\ali{sl}{2}$\fi}

\renewcommand{\Im}{\mathop{\rm Im}\nolimits}

\newcommand{\pmat}[1]{\begin{pmatrix}#1\end{pmatrix}}
\newcommand{\smat}[1]{\left(\begin{smallmatrix}#1\end{smallmatrix}\right)}

\newcommand{\bMA}[1]{\[\begin{array}{#1}}
\newcommand{\eMA}{\end{array}\]}
\newcommand{\HF}{\hfill}
\newcommand{\ol}[1]{\overline{#1}}
\newcommand{\C}{\NC}
\newcommand{\BC}[1]{\NC^{#1}}

\newcommand{\NR}{{{\mathbb R}}}

\newcommand{\NC}{{{\mathbb C}}}

\newcommand{\NZ}{{{\mathbb Z}}}

\def\be{\begin{equation}}
\def\ee{\end{equation}}

\def\C{{\mathbb C}}

\def\fd{{f^\dagger}}

\def\cpn{{\C P^{N}}}
\def\bp{{\bar{\partial}}}
\def\p{{\partial}}

\def\d{{\mathrm{d}}}

\def\bw{{\bar{W}}}
\def\tr{{\mathrm{tr}}}

\mathchardef\za="710B  
\mathchardef\zb="710C  
\mathchardef\zg="710D  
\mathchardef\zd="710E  
\mathchardef\ze="710F  
\mathchardef\zz="7110  
\mathchardef\zh="7111  
\mathchardef\zy="7112  
\mathchardef\zi="7113  
\mathchardef\zk="7114  
\mathchardef\zl="7115  
\mathchardef\zm="7116  
\mathchardef\zn="7117  
\mathchardef\zx="7118  
\mathchardef\zp="7119  
\mathchardef\zr="711A  
\mathchardef\zs="711B  
\mathchardef\zt="711C  
\mathchardef\zu="711D  
\mathchardef\zf="711E  
\mathchardef\zq="711F  
\mathchardef\zc="7120  
\mathchardef\zw="7121  
\mathchardef\zG="7000  
\mathchardef\zD="7001  
\mathchardef\zY="7002  
\mathchardef\zL="7003  
\mathchardef\zX="7004  
\mathchardef\zP="7005  
\mathchardef\zS="7006  
\mathchardef\zU="7007  
\mathchardef\zF="7008  
\mathchardef\zC="7009  
\mathchardef\zW="700A  
\mathchardef\ze="7122  
\mathchardef\zvy="7123  
\mathchardef\zvr="7125 
\mathchardef\zvs="7126 
\mathchardef\zvf="7127  

\begin{document}


%
\renewcommand{\evenhead}{A~M Grundland, A~Strasburger and W~J~Zakrzewski}
\renewcommand{\oddhead}{Surfaces immersed in $\su{N+1}$ Lie algebras}

%
\thispagestyle{empty}

\FirstPageHead{*}{*}{2005}{\pageref{firstpage}--\pageref{lastpage}}{Article}

\copyrightnote{2005}{A~M Grundland, A~Strasburger and W~J~Zakrzewski}

\Name{Surfaces immersed in $\su{N+1}$ Lie algebras obtained from the $\cpn$ sigma
models}

\label{firstpage}

\Author{A~M GRUNDLAND~$^{\dag^1}$, A~STRASBURGER~$^{\dag^2}$ and W~J~ZAKRZEWSKI~$^{\dag^3}$}

\Address{$^{\dag^1}$ Centre de Recherches Math{\'e}matiques, Universit{\'e} de Montr{\'e}al, \\
C. P. 6128, Succ.\ Centre-ville, Montr{\'e}al, (QC) H3C 3J7, Canada \\
Universit\'{e} du Qu\'{e}bec, Trois-Rivi\`{e}res CP500 (QC) G9A 5H7, Canada \\
~~E-mail: grundlan@crm.umontreal.ca\\[10pt]
$^{\dag^2}$ Department of Mathematical Economics, Warsaw Agricultural University,\\
ul. Nowoursynowska 166, 02-787 Warszawa, Poland \\
~~E-mail: strasburger@alpha.sggw.waw.pl\\[10pt]
$^{\dag^3}$ Department of Mathematical Sciences, University of Durham,\\
Durham DH1 3LE, United Kingdom \\
~~E-mail: W.J.Zakrzewski@durham.ac.uk}

\Date{Received Month *, 2005; Revised Month *, 2005; 
Accepted Month *, 2005}



\begin{abstract}
We study some geometrical aspects of two dimensional orientable surfaces arrising from the study of $\cpn$ sigma models. 
To this aim we employ an identification of $\NR^{N(N+2)}$ with the Lie algebra  
$\su{N+1}$ by means of which we construct a generalized Weierstrass formula for immersion of such surfaces.  The structural elements of the surface like its moving frame, the Gauss-Weingarten and the Gauss--Codazzi--Ricci
equations are expressed in terms of the solution of the $\cpn$ model defining it. 
Further, the first and second fundamental forms, the Gaussian 
curvature, the mean curvature vector, the Willmore functional and
the topological charge of surfaces are expressed in terms of this solution. 
We present detailed implementation of these results for surfaces immersed in $\su{2}$ and $\su{3}$ Lie algebras.
\end{abstract}

\medskip

\section{Introduction}
\label{sec:Intro}
In this paper we develop further the study of immersions of two-dimensional
surfaces in multidimensional Euclidean or Minkowski spaces by means of $\cpn$ models,
which was carried out in a series of papers \cite{GZ01,GRZA,GZ,GS1}.
The key point is the formulation of the equations defining the immersion
directly in the matrix form (cf. equation \ref{3_2}) where the immersion
takes values in the Lie algebra $\su{N+1}$, identified by means of the negative of the Killing form with the Euclidean space $\NR^{N(N+2)}$. This allows us to formulate explicitly the structural equations for the immersion (the Gauss--Weingarten 
and the Gauss--Codazzi--Ricci equations) directly in matrix terms. In particular in 
Proposition \ref{prop3}  we establish an explicit form of the Gauss-Weingarten equations
satisfied by the moving frame on a surface corresponding to the $\cpn$ model.
This is done in a fashion independent of any specific parametrization. Then we use this result to  establish various geometric characteristics of the studied immersions like curvatures and curvature vectors. All these quantities are directly derived from the map describing the relevant $\cpn$ model. 
For the simplest case $N=1$ the equation \Ref{3_2} defining the immersion takes
the form 
\begin{equation}
dX = i(dX_1 \zs_2 + dX_2 \zs_1 + dX_3 \zs_3),
\end{equation}
where $\zs_1,\zs_2,\zs_3$ are the usual Pauli matrices, and the differentials
of coordinate functions of the immersion $dX_1,dX_2,dX_3$ are given in terms 
of the affine coordinate $W$ of ${\C P^{1}}$ model by equation \Ref{4_3}.
As follows from the equation \Ref{1_22'} describing the ${\C P^{1}}$ model, a particular class of solutions of this model is given by an arbitrary holomorphic function $W$ --- in this case the immersion is minimal (i.e. represents a surface with zero mean curvature) and $W$ expresses the Gauss map of the surface by means of the stereographic projection. This is directly related to the classical Weierstrass--Enneper formulae for an immersion of a minimal surface in $\NR^3$. In fact, almost one and half century ago Weierstrass
and Enneper showed \cite{Weier,Enne}  that every minimal surface in $\NR^3$ can be represented locally in terms of two holomorphic functions $\psi$ and $\phi$ defined on a domain $\mathbb{D}\in \C$ by the following expressions 
\begin{equation}
X(\xi,\bar{\xi}) = \mathrm{Re} \left( \int_0^{\xi} (\psi^2-\phi^2) \d \xi' , i\,
\int_0^{\xi} (\psi^2+\phi^2) \d \xi' , -2 \int_0^{\xi} \psi \phi
\d \xi' \right).
\label{minimal}
\end{equation}
This implies that the complex tangent vector of the immersion is given by
\begin{equation}
\p X_1 = \psi^2-\phi^2,\quad \p X_2 = i\,(\psi^2+\phi^2),\quad \p X_3 =
-2\psi \phi, 
\end{equation}
where $\p$ denotes the (complex) derivative with respect to $\xi$. Moreover, the metric of the minimal surface is conformal and is expressed in terms of local parameters $\xi$ and $\bar{\xi}$  by the formula 
\begin{equation}
\d s^2 = 2(|\psi|^2+|\phi|^2)^2 \d \xi \d \bar{\xi}. 
\end{equation}
This implies, in particular, that the coordinate lines $\xi=\text{const}$ and $\bar{\xi}=\text{const}$ describe geodesics on this surface.  

The ideas of Enneper and Weierstrass 
have since been developed by many authors with a purpose of extending them to construct immersions of more general types of surfaces.  For a classical presentation we refer to a treatise by Eisenhardt\cite{Ei} and for a modern approach to the subject see e.g. \cite{Bob,BobE,Fok,FokG,Ken,Osser,Ped,Pink1,Pink2} and references therein, in particular the recent books by F. Helein \cite{Hel,HelF} and K. Kenmotsu \cite{Ken2}.  This topic has been further explored by, among others, B. Konopelchenko \cite{Kono1}{}. He established a direct connection between generic surfaces and trajectories of an infinite-dimensional Hamiltonian system.  Namely, he considered a
nonlinear Dirac type system of equations for two complex-valued
functions $\psi_1$ and $\psi_2$ of $\xi, \bar{\xi}$
\begin{equation}
\label{0_4}
\p \psi_1 = (|\psi_1|^2+|\psi_2|^2)\psi_2, \qquad 
\bp \psi_2 = -(|\psi_1|^2+|\psi_2|^2)\psi_1.
\end{equation}
He showed \cite{Kono2} that for any solutions $\psi_1,\,\psi_2 $ of the system \Ref{0_4} and an arbitrary curve $\gamma$ in the complex plane, the following integrals over the bilinear combinations of $\psi_i$, $i=1,2$  
\begin{equation}
\label{0_5}
\begin{split}
X_1 &= \int_{\gamma} \, (\psi_1^2 - \psi_2^2) \d \xi' +
(\bar{\psi}_1^2 - \bar{\psi_2}^2) \d \bar{\xi}',\\
X_2 &= \int_{\gamma} (\psi_1^2+\psi_2^2) \d \xi' - (\bar{\psi}_1^2
+ \bar{\psi_2}^2) \d \bar{\xi}',\\
X_3 &= -\int_{\gamma} \psi_1\psi_2 \d \xi +
\bar{\psi}_1\bar{\psi}_2 \d \bar{\xi}',
\end{split}
\end{equation}
determine the coordinates of the radius vector $X=(X_1,\,X_2,\,X_3)$ describing 
a constant mean curvature (CMC) surface immersed in $\mathbb{R}^3$ (the
integrals are independent of the integration path, since the integrands are
exact differentials). To see how to reduce more general cases down to this case see e.g. Ref \cite{KL,Kono2}{}.  In accordance with \cite{Carr} we will refer to equations \Ref{0_4}
and \Ref{0_5} as the generalized Weierstrass formulae. 

It was shown in our previous work \cite{GZ} that the generalized 
Weierstrass formulae for two-dimensional surfaces with non-vanishing mean curvature in multi-dimensional spaces are equivalent to  $\cpn$ sigma models. This determination has opened a new way for constructing and studying two-dimensional surfaces.  The further advantage of use 
of the $\cpn$ models in this context lies in the fact that they
allow us to replace the methods based on Dirac-type equations by
the formalism connected with completely integrable systems, for
example Lax pairs, Hamiltonian structures, or systems defining
infinite number of conserved quantities.  An original procedure for
constructing the general classical solutions admitting finite
action of the Euclidean two-dimensional  $\cpn$  model was devised
by A.M. Din et al \cite{Din}.  These solutions are obtained by repeated applications of a certain transformation to the basic solution expressed in terms of holomorphic
functions. As a result, one gets three classes of solutions:
holomorphic, anti-holomorphic and the ``mixed'' ones. In this paper we show that to each of these  solutions we can associate a surface in $\su{N+1}\simeq \NR^{N(N+2)}$. In the holomorphic (or antiholomorphic) case we are able to integrate completely the equations of  the immersion. It turns out that in the $\C P^1$ case the surface is a part of an Euclidean sphere, cf. Section \ref{sec:CP1}.  However for arbitrary $N>1$ other situations are possible.  In Example $1$ in Section \ref{sec:CP2} we present a one-parameter family of surfaces for which the curvature is not constant but for some specific values of this parameter it reduces to a constant. 

The second and third examples in Section \ref{sec:CP2} are concerned with mixed solutions of the $\C P^2$ model. In one case we obtain a surface in $\NR^8$, which happens to be immersible in $\NR^3$, but the immersion does not come from a $\C P^1$ model, since the curvature is not constant.  The other mixed solution leads to a generic surface in $\NR^8$ with nonconstant curvatures. All these results raise interesting questions which require further investigations concerning general properties of immersions given by the $\cpn$ models - either holomorphic or mixed. 

Finally, let us note that the outlined approach to the study of 
surfaces in $\NR^m$ lends itself to numerous potential applications. It is useful for  description of monodromy of solutions of higher order Painlev{\'e} equations and their connection with theory of surfaces.  It can also lead to the development of numerical computing tools in the study of surfaces through the techniques of completely integrable systems.  

Surfaces immersed in Lie groups, Lie algebras and homogeneous
spaces appear in many areas of physics, chemistry and biology\cite{PolS,Nel,David,SomF,Chen,Dav,Ou,Saf,Lan}{}.  The algebraic approach to structural equations of these surfaces has often proved to be very difficult from the computational point of view.
A natural geometric approach to derivation and classification of
such equations which we propose here seems therefore to be of importance for
applications in physics and other sciences. 

This paper is organized as follows. In Section \ref{sec:CPN} we introduce basic material on $\cpn$ models --- in presenting it we focus on the use of a compatibility condition, rather then on the usual approach via the Euler--Lagrange equations. The next section is devoted to the  presentation of required notions and facts on the structure of complex projective spaces. Here we prove in detail a certain decomposition of the group $\gr{SU}{N+1}$, which was previously noted in the paper \cite{Rowe} of Rowe et al.  In Section \ref{sec:SUN} we show how to use the equations of the $\cpn$ model to construct an immersion of a surface in the Lie algebra $\su{N+1}$. The obtained formula extends the classical Weierstrass formula for the conformal immersion of a $2$-dimensional surface into the $3$-dimensional Euclidean space. We derive the equations of the moving frame for the above immersion in $\su{N+1}$ and derive some geometrical characteristics of these surfaces. This analysis is developed further in Sections \ref{sec:CP1} and \ref{sec:su3}, where, using the conformality of the surfaces obtained from the $\cpn$ models for $N=1,\,2$, we establish an explict formula  for the moving frame in terms of the data of the model. Let us note that the case $N=2$ produces surfaces immersed in $\NR^8\simeq\su{3}$. In Section \ref{sec:CP2} we illustrate our theoretical considerations with some examples based on explicit solutions of the $\C P^2$ sigma model.  The last section contains remarks and sugestions regarding possible further developments.  

\section{Preliminaries on $\cpn$ models \cite{WZ}}
\label{sec:CPN}
From a large supply of geometrical models of immersions \cite{GU1,GU2} we  
concentrate in this paper on a particular class of models, the so
called $\cpn$ sigma models.  The $\cpn$ sigma model can be defined in
terms of functions
\begin{equation}
\C \supset \Omega \ni \xi = \xi^1+i\,\xi^2 \mapsto 
z=(z_0,\,z_1,\,\ldots,\,z_N) \in \C^{N+1} 
\end{equation}
defined on an simply connected domain (i.e. an open connected subset) $\Omega$ of the complex plane and satisfying the constraint $z^{\dagger}\cdot z=1$. Here and below we employ the standard notation where points of the complex coordinate space $\BC {N+1}$ are denoted by $z=(z_0,\, z_1,\,\ldots,\,z_N)$ and the standard hermitian inner product in $\BC{N+1}$, respectively the norm, by 
\begin{equation} 
z^\dagger\cdot w = \langle{z,\,w}\rangle =\IS{j=0}{N}\overline{z}_jw_j,
\qquad \text{resp.}\quad |z|=(z^\dagger\cdot z)^{1/2}.
\end{equation}
The unit sphere in $\BC{N+1}$ corresponding to this norm is given by 
\[
S^{2N+1}=\{z\in \BC{N+1}\mid z^\dagger\cdot z= |z|^2=1\}.
\]
We shall use $\partial_{\mu} = \partial/\partial \xi^{\mu},\ \mu=1,2$ to denote ordinary derivatives and $D_{\mu}$ for the covariant derivatives defined according to the formula 
\begin{equation}
D_{\mu} z = \partial_{\mu} z - z(z^{\dagger}\cdot \partial_{\mu} z).
\end{equation} 
With this notation the Lagrangian density for such a model is given by (cf. e.g. \cite{WZ})
\begin{equation}
\label{1_77}
\mathcal{L}=\frac{1}{4}\left(D_{\mu}z\right)^{\dagger} \cdot
(D_{\mu}z),
\end{equation}
and the solutions of the $\cpn$ model are stationary points of the action
functional
\begin{equation}
\label{1_4} S=\int_{\Omega}\mathcal{L}\,\d\xi\d\bar{\xi} =\frac{1}{4} \int_{\Omega}
\left(D_{\mu}z\right)^{\dagger} \cdot (D_{\mu}z)\,\d\xi\d\bar{\xi}.
\end{equation}
The physically relevant case concerns fields which can be extended to the whole
Riemann sphere $S^2=\C \cup \{\infty\}$, but the case of an arbitrary $\Omega$ is also of some interest, especially for questions of (local) differential geometry.

Next we note that since $\mathcal{L}$ is not changed by the transformations $z\mapsto e^{i\varphi}z$ with $\varphi\in \NR$, it is actually defined by the map $[z] : \Omega \rightarrow {\C P^{N}}$, where  $[z]= \{e^{i\varphi}z\mid \varphi \in \NR\}$ is the element of the projective space ${\C P^{N}}$ corresponding to $z\in S^{2N+1}$. We find it often more convenient to use this latter point of view and describe the model in terms of ``unnormalized" fields  $ \xi\mapsto f = (f_0,\,f_1,\ldots,\,f_N) \in \C^{N+1}\setminus\{0\}$ related to the ``$z$'s'' above by 
\begin{equation}
z=\frac{f}{|f|}, \quad \text{where}\quad  |f|=(\fd \cdot f)^{1/2}.
\end{equation}
We shall refer to the ``$z$'s'' above as inhomogeneous and to the ``$f$'s'' as homogeneous coordinates of the model. 

Now, using the customary notation for holomorphic and
antiholomorphic derivatives
\begin{equation}
\partial = \frac{1}{2}\left(\frac{\partial}{\partial
\xi^1}-i\frac{\partial}{\partial \xi^2}\right) \quad , \quad
\bar{\partial} = \frac{1}{2}\left(\frac{\partial}{\partial
\xi^1}+i\frac{\partial}{\partial \xi^2}\right)
\end{equation}
and introducing the orthogonal projector on the orthogonal complement to the complex line in $\mathbb{C}^{N+1}$ determined by $f$ given by   
\begin{equation}
P=\ID_{N+1} -\frac{1}{\fd\cdot f} f \otimes \fd
\end{equation}
we may express the action functional \Ref{1_4} in terms of $f$'s  by 
\begin{equation}
\label{1_7} S=\frac{1}{4}\int_{\Omega} \frac{1}{\fd\cdot f}(\partial
\fd P \bar{\partial}f + \bar{\partial} \fd P \partial f) \d\xi
\d\bar{\xi},
\end{equation}
Since $P$ is an orthogonal projector, it satisfies 
\begin{equation}
\label{1_9} P^2=P, \quad \quad P^{\dagger} = P.
\end{equation}
The map $[z]$ is determined by a solution of the Euler--Lagrange equations which is associated with the action \Ref{1_7}.  In terms of homogeneous coordinates $f$'s the equations take the form
\begin{equation}
\label{1_11}
P\left[\partial \bar{\partial f} - \frac{1}{\fd\cdot f}
\left(\left(\fd \cdot \bar{\partial} f\right)\partial f +
\left(\fd \cdot \partial
f\right)\bar{\partial} f\right)\right] = 0. 
\end{equation}
Using the projector P we can rewrite \Ref{1_11} as 
\begin{equation}
\label{1_12}
 \left[\partial\bar{\partial}P,P\right]=0,
\end{equation}
or equivalently as the conservation law
\begin{equation}
\label{1_13}
\partial
\left[\bar{\partial}P,\,P\right]+\bar{\partial}\left[\partial P,\,P\right] = 0.
\end{equation}
Further, introducing the $(N+1)$ by $(N+1)$ matrix $\mathbb{K}$  
\begin{equation}
\label{1_14} \mathbb{K}=\left[\bar{\partial}P,P\right] =
\frac{\bar{\partial} f \otimes \fd - f\otimes
\bar{\partial}\fd}{\fd\cdot f} + \frac{f\otimes \fd}{(\fd \cdot f)^2}
\left[(\bar{\partial}\fd \cdot f) - (\fd \cdot \bar{\partial}
f)\right]
\end{equation}
and noting that its hermitian conjugate is 
\begin{equation}
\mathbb{K}^{\dagger}= -\left[\partial P,P\right] = \frac{-\partial f
\otimes \fd + f\otimes \partial \fd}{\fd\cdot f} + \frac{f\otimes
\fd}{(\fd\cdot f)^2} \left[(\partial \fd \cdot f) - (\fd \cdot
\partial f)\right], 
\end{equation}
we can reformulate \Ref{1_13} succinctly as 
\begin{equation}
\label{1_16}
\partial \mathbb{K} - \bp \mathbb{K}^{\dagger}=0.
\end{equation}
It follows from the above \Ref{1_16} that $\partial \mathbb{K}$ is an hermitian matrix, i.e. $\partial \mathbb{K} \in i\su{N+1}$.  

One can check by a straightforward computation that the
complex-valued functions
\begin{equation}
\label{1_18}
J=\frac{1}{\fd\cdot f} \partial \fd P \partial f, \qquad  
\bar{J} = \frac{1}{\fd \cdot f} \bar{\partial} \fd P \bar{\partial} f,
\end{equation}
satisfy
\begin{equation}\label{1_26}
\bar{\partial}J = 0 \quad , \quad \partial\bar{J}=0
\end{equation}
for any solution $f$ of the Euler--Lagrange equations \Ref{1_11}.
Note that $J$ and $\bar{J}$ are invariant under
global $\gr{U}{N+1}$ transformation, i.e. $f \rightarrow \psi f,\, \psi
\in \gr{U}{N+1}$.\\

\section{Some decompositions of $\gr{SU}{N+1}$ and related parametrizations of $\cpn$}
\label{sec:SUN}
In this section we collect several facts concerning realization of projective space 
$\cpn$ as a homogeneous space of the special unitary group $\gr{SU}{N+1}$ and discuss related decompositions of the group and its Lie algebra $\su{N+1}$. 
The standard reference, where all the details missing here can be found, is the book of Helgason\cite{Helg}{}. 

As is well known, the space $\cpn$  consists of complex lines through the origin $0$  in $\C^{N+1}$ (the one dimensional complex subspaces of $\C^{N+1}$). We denote by $\pi$ the map which associates to any nonzero vector $Z=(z_0,\,\ldots,\,z_N)\in \C^{N+1}$ the line passing through the origin and $Z$,  so that 
\begin{equation}
\pi(Z)=\{\zl Z\mid \zl\in \NC\}=[z_0,\,\LI{N}{z}],    
\end{equation}
The numbers $z_0,\,\LI{N}{z}$, determined up to a nonzero complex number, are called the homogeneous coordinates of the line $\pi(Z)$. The restriction of $\pi$ to the unit sphere $S^{2N+1}=\{Z\in\BC{N+1}\mid Z^\dagger Z = 1\}$ 
remains surjective --- the resulting map $\pi_H :S^{2N+1}\to \NC{P^N}$ is known as the Hopf fibration. Observe that if the line $l$ passes through the point $Z_0\in S^{2N+1}$, then the fiber over $l$, $\pi_H^{-1}(l)=\{Z\in S^{2N+1}|Z= e^{i\zvf}Z_0, \zvf\in \NR\}$, is just the great circle in $S^{2N+1}$ passing through $Z_0$. 

For any given $j = 0,\,1,\,\ldots,\,N$ one introduces the so called affine or inhomogeneous coordinates  defined in the complement of the set $H_j=\{\pi(Z)\mid Z\in \NC^{N+1}_*,\  z_{j}=0\}\subset \NC{P^N}$ by the prescription 
$$
[z_0,\,\LI{N}{z}]\mapsto \biggl(\frac{z_{0}}{z_{j}},\,\ldots,\,\frac{z_{j-1}}{z_{j}},\,\frac{z_{j+1}}{z_{j}},\,\ldots,\,\frac{z_{N}}{z_{j}}\biggr)
$$ 
which sets up a natural isomorphism of $\NC^N$ with $\NC{P^N}\setminus H_j$. In the particular case of the affine coordinates 
defined in the set $U_0 =  \NC{P^N}\setminus H_0$ we shall write $W_i= z_{i}/z_{0}$.   

By transitivity of the action of $\gr{SU}{N+1}$ on the set of lines in $\C^{N+1}$  one has a natural identification $G/K_0 \simeq \NC{P^N}$, with $K_0$ denoting the isotropy group of the standard reference point $l_0=\pi_H(e_{0})=\NC e_{0}$. 
Now  
\begin{equation}
K_0  =\mathbf{S}(\gr{U}{1}\times \gr{U}{N})= 
\left\{\pmat{(\det{\mathbf{A}})^{-1} & 0\\  0 & \mathbf{A}}
\mid \mathbf{A}\in \gr{U}{N}\right\}
\end{equation}
and the identification above is written  as
$\NC{P^N} \simeq \gr{SU}{N+1}/ \mathbf{S}(\gr{U}{1}\times \gr{U}{N})$.
Passing to the Lie algebra level and denoting the respective Lie algebras by $\mathfrak{g}=\ali{su}{N+1}$ and 
\[
\mathfrak{k}_0=
\left\{\pmat{-\Tr \mathbf{A} & 0\\ 0 & \mathbf{A} }\mid \mathbf{A}\in \ali{u}{N}
\right\}
\] 
one has the direct sum decomposition of the isotropy Lie algebra 
\[
\mathfrak{k}_0 =\mathfrak{c} \oplus  \ali{su}{N},\qquad \text{where}\quad 
\mathfrak{c}= 
\left\{\pmat {-\zm &0\\ 0 & \dfrac{\zm}{N}\ID_{N}}\mid \zm\in i\NR\right\}\simeq i\NR
\] 
with $\mathfrak{c}$ being its center. To study other decompositions we first recall  
that the Killing form of $\mathfrak{g}=\ali{su}{N+1}$  is given by the formula \begin{equation}\label{eqn:Killing}
  B(X,\,Y)=2(N+1)\Tr(XY)
\end{equation}
and is negative definite. The space $\ali{su}{N+1}$ of skew-hermitian matrices can thus be given the structure of a real Euclidean space of dimension $N(N+2)$ by taking the negative of the Killing form as the inner product.   
The orthogonal complement to $\mathfrak{k}_0$ with respect to this inner product consists of matrices of the form  
\begin{equation}\label{eqn:Z(x)}
Z(x) = \pmat{0 &  -x^\dagger \\ x & \mathbf0_N } = 
x\otimes e_{0}^\dagger - e_{0}\otimes x^\dagger,\qquad 
\end{equation}
where $x =(x_1,\,\ldots,\,x_N)\in \NC^N$ and $\mathbf0_N$ is the $N\times N$ zero matrix, and this yields the orthogonal decomposition 
$\mathfrak{g}=\mathfrak{k}_0\oplus\mathfrak{p}$.  

This later fact is a starting point for introducing a useful parametrization of the projective space, analogous to the spherical coordinates on the Euclidean sphere.  Observe that the adjoint action of the isotropy group $K_0$ on $\mathfrak{p}$ reduces to the action of $\gr{U}{N}$ on $\NC^N$ given by the following formula 
\[
\pmat{(\det{A})^{-1} & 0\\ 0 & A}
\pmat{0 &  -x^\dagger \\ x & \mathbf0_N }\pmat{\det{A} & 0\\ 0 & A^{-1}}= 
\pmat{0 &  -(\det{A})^{-1}(Ax)^\dagger \\ \det(A)Ax & \mathbf0_N }
\]
The action is clearly transitive on the unit sphere in $\mathfrak{p}$ and essentially this fact implies validity of the next result (for a general form of such decompositions cf. \cite[p.~402]{Helg}).
To formulate it we first introduce more notations.
Set $H = Z(e_1)$ 
and let $\mathfrak{a}=\NR{H}\subset \mathfrak{p}$. The Lie subgroup $\exp\mathfrak{a} =\{\exp\za H\mid \za\in \NR\}$ of $\gr{SU}{N+1}$ consists of matrices 
\begin{equation}\label{eqn:max_torus1}
\exp\za H= \pmat{\cos\za & -   \sin\za & 0 \HF \\
   \sin\za & \HF \cos\za & 0 \HF \\ 
                 0       &  0          & \ID_{N-1} } 
                 = \pmat{\mathbf{R}(za) &0 \\0&  \ID_{N-1}},  
\end{equation}
where we have set $\mathbf{R}(\za)= \smat{\cos\za & - \sin\za \\ \sin\za & \HF \cos\za}$ and which is isomorphic to $\gr{SO}{2}$. The corresponding maximal torus in $\cpn$ is 
\begin{equation}\label{eqn:max_torus2}
A = (\exp{\mathfrak{a}})\cdot l_0 = \left \{\bigl[\cos\za,\,\sin\za,\,0,\, \ldots,\, 0 \bigr] \mid \za\in \NR\right \}
\end{equation}

Denoting further by $M\subset K_0$ the centralizer and by $M'\subset K_0$ the normalizer of $\mathfrak{a}$ in $K_0$, i.e. $M=\{k\in K_0\mid kHk^{-1}=H\}$ and 
$M' =\{k\in K_0\mid kHk^{-1}\subset \NR H \}$, we see that  
\begin{gather}
M= \left\{\pmat {u & 0& 0 \\ 0 & u & 0 \\ 0 & 0 &U } 
\Biggm | U\in \gr{U}{N-1},\, u\in \gr{U}{1},\ u^2\det{U}=1 \right\} \\
M' = \left\{\pmat {u & 0 & 0 \\ 0 & \ze{u} & 0 \\ 0 & 0 & \ze{U} } 
\Biggm | U\in \gr{U}{N-1},\, u\in \gr{U}{1}, \ \ze=\pm1,\ u^2\det{U}=1\right\}
\end{gather}
The factor group $M'/M =\NZ_2$ is the Weyl group associated with $\cpn$.  

We can now formulate the main result of this section, which will be used extensively later on. It describes a certain decomposition of the group $\gr{SU}{N+1}$ related to the spherical parametrization of the projective space and is stated as the point {\bf c)} of the proposition below. The points {\bf a)} and {\bf b)} are included for readers' convenience and comprise the classical decompositions found e.g. in  \cite[p.~402]{Helg}). It should be pointed out that {\bf c)}  is an elaboration of the result stated in \cite{Rowe}{}, but proved there only for $\gr{SU}{3}$. 

\begin{theo}[Polar decompositions]\label{theo:Polar decomp.}\hfill\\ 
Let $G=\gr{SU}{N+1}$, $K_0=\mathbf{S}(\gr{U}{1}\times \gr{U}{N})$, and let $K_1$ denote the image of $\gr{SU}{N}$ in $K_0$ by means of the injection $\mathbf{A}\mapsto \smat{1&0\\ 0&\mathbf{A}}$. The remaining notations are as explained above.   

\noindent {\bf a)} Every element of $G$ can be written as a product $g=k_1(\exp{\zy{H}})k_2$ with $k_1,\,k_2\in K_0$. More precisely, the map 
\begin{equation}\label{eqn:Cartan_decomp}
K_0\times \exp{\mathfrak{a}}\times K_0 \ni (k_1,\,\exp{\zy{H}},\,k_2)
\mapsto k_1(\exp{\zy{H}}){k_2}\in G
\end{equation} 
arising from the group multiplication, is a smooth surjection. 

\noindent {\bf b)} The map 
\begin{equation}\label{eqn:polar_decomp}
K_0/M \times A \ni (kM,\,\exp{\zy{H}}\cdot l_0)\mapsto k\cdot \exp{\zy{H}}\cdot l_0\in \NC{P^N}  
\end{equation}
is a smooth surjection and is a double covering on the complement of the point $l_0\in\NC{P^N}$.  

\noindent {\bf c)} Let 
\begin{equation}\label{eqn:delta}
\zd: \gr{U}{1}\ni \zm \mapsto \zd(\zm,\,\ol{\zm}) =
\smat{\zm& 0\\0& \overline{\zm}}
\end{equation}
be the diagonal embedding of $\gr{U}{1}$ into $\gr{SU}{2}$. 
The map 
\begin{multline}\label{eqn:Rowe_decomp} 
\gr{SU}{N} \times \gr{U}{1}\times \gr{SO}{2} \times \gr{SU}{N}\quad \longrightarrow \quad \gr{SU}{N+1} \\
(\mathbf{A}_1,\,\zm,\,\mathbf{R}(\zy),\,\mathbf{A}_2) \mapsto \\
\pmat{1&0\\0& \mathbf{A}_1}\cdot
\pmat{\zd(\zm,\,\ol{\zm}) & 0 \\ 0 &\ID_{N-1}} \cdot
\pmat{\mathbf{R}(\zy) & 0\\ 0&\ID_{N-1}} \cdot
\pmat{\zd(\zm,\,\ol{\zm}) & 0 \\ 0 &\ID_{N-1}} \cdot
\pmat {1 & 0\\ 0 & \mathbf{A}_2}
\end{multline}
is a smooth surjection. 
\end{theo}
\begin{proof}
We are going to prove only part {\bf c)} of the statement and this follows by simple matrix calculations from the polar decomposition given in equation~\Ref{eqn:Cartan_decomp}. Assume that we have a product of the form
\begin{equation}\label{eqn:prelim_prod}
  \pmat{\zl_1&0\\0& \mathbf{A}_1}\cdot
\pmat{\mathbf{R}(\zy) & 0\\ 0&\ID_{N-1}} \cdot
\pmat {\zl_2 & 0\\ 0 & \mathbf{A}_2},
\end{equation}
where $\mathbf{A}_i\in \gr{U}{N}$ and $\zl_i\det\mathbf{A}_i=1$ for $i=1,\,2$. By splitting a factor of the form 
\[
\pmat{\zd_3(\zl,\,\zl,\,\zl^{-2}) & 0 \\ 0 &\ID_{N-2}},\qquad \zl=\zl_1^{1/2}\zl_2^{-1/2} 
\]
from the matrix on the left hand side in this product and commuting it with the middle term, we can bring the entries in the top left corners of the matrices on both sides of  equation~\Ref{eqn:prelim_prod} to coincide with each other, thus obtaining  the product 
\[
\pmat{\zl&0\\0& \mathbf{A'}_1}\cdot
\pmat{\mathbf{R}(\zy) & 0\\ 0&\ID_{N-1}} \cdot
\pmat {\zl & 0\\ 0 & \mathbf{A'}_2}
\]
with the relations $\zl\det\mathbf{A'}_i=1$ for $i=1,\,2$ still satisfied. Now it remains only to split off the factors of the form $\smat{\zd(\zl,\,\ol{\zl}) & 0 \\ 0 &\ID_{N-1}}$ from the both extreme terms to get the sought for form \Ref{eqn:Rowe_decomp} of the product. 
\end{proof}
\begin{rem}\label{rem:product}\normalfont
The statement in {\bf c)} above is simply that every element of $\gr{SU}{N+1}$ can be written as a product of four matrices of the above given form, and 
writing down the product of the middle terms in equation~\Ref{eqn:Rowe_decomp} explicitely we obtain  
\begin{multline}
\pmat{\zd(\zm,\,\ol{\zm}) & 0 \\ 0 &\ID_{N-1}} \cdot
\pmat{\mathbf{R}(\zy) & 0\\ 0&\ID_{N-1}} \cdot
\pmat{\zd(\zm,\,\ol{\zm}) & 0 \\ 0 &\ID_{N-1}} = \\
\noalign{\medskip}
\pmat{\zm^2\cos\zy & -   \sin\zy & 0 &  \ldots & 0 \\ 
       \sin\zy & \HF \zm^{-2}\cos\zy & 0 &  \ldots & 0 \\
   \noalign{\smallskip}
  {\displaystyle0\atop\displaystyle\vdots}&{\displaystyle0\atop\displaystyle\vdots}&& 
 \ID_{N-1}& \\
  0  & 0  & & &} 
\end{multline}
so that the right hand side of the decomposition \Ref{eqn:Rowe_decomp} reduces to 
\begin{equation}
\pmat{ 1 & 0 &  \ldots & 0\\ 
\noalign{\smallskip}
{\displaystyle 0 \atop \displaystyle  \vdots} & & \mathbf{A}_1  &\\
\noalign{\smallskip}
0  & &   & }
\pmat{\zm^2\cos\zy & -   \sin\zy & 0 &  \ldots & 0 \\
       \sin\zy & \HF \zm^{-2}\cos\zy & 0 &  \ldots & 0 \\ 
   \noalign{\smallskip}
  {\displaystyle0\atop\displaystyle\vdots}&{\displaystyle0\atop\displaystyle\vdots}&& 
 \ID_{N-1}& \\
  0  & 0  & & &} 
\pmat{ 1 & 0 & \ldots & 0\\ 
{\displaystyle 0 \atop \displaystyle \vdots}& & \mathbf{A}_2 &\\
\noalign{\smallskip}
0 &   & },
\end{equation}
where $\mathbf{A}_1,\,\mathbf{A}_2\in \gr{SU}{N},\ \zm=e^{i\za } \in \gr{U}{1}$ and $\zy,\,\za\in \NR$, which up to unimportant changes in parametrization, is the expression given in \cite[Equation~(2) on p.~3605]{Rowe}.
\end{rem}

The polar decomposition \Ref{eqn:Cartan_decomp} reduces to the following decomposition of $\gr{SU}{3}$ which can be found in the paper of D.~J.~Rowe and all \cite{Rowe}.
\begin{coro}\label{coro:Cartan_for_SU3} Each element of the $\gr{SU}{3}$ group can be decomposed into the following product:
\begin{equation}
g= \pmat{1& 0 \\ 0& U_1 }\pmat{\zd(\zl,\,\ol{\zl})&0\\0& 1 }
\pmat{R(\za)&0\\0& 1 }\pmat{\zd(\zl,\,\ol{\zl})&0\\0& 1 }\pmat{1& 0 \\ 0& U_2},
\end{equation}
where $U_i\in \gr{SU}{2}$ for $i=1,\,2$ and 
$$ R(\za)=\pmat{\cos{\za}& -\sin{\za}\\ \sin{\za}& \cos{\za} }\quad  \text{and}\quad  
\zd(\zl)=\pmat{\zl& 0\\0& \overline{\zl}}
$$
with $\zl\in \NC$ with $|\zl|=1$.
Writing this more explicitely we obtain 
\begin{equation}
g= \pmat{1&0&0 \\0 &\HF a_1 & b_1 \\0& -\overline{b_1} & \overline{a_1} }
\pmat{\zl & 0 & 0 \\ 0&\overline{\zl}&0\\0 &0 & 1 }
\pmat{\cos{\za}& -\sin{\za} &0 \\ \sin{\za}&\HF \cos{\za}&0 \\ 0 &0 & 1}
\pmat{\zl & 0 & 0 \\ 0&\overline{\zl}&0\\0 &0 & 1 }
\pmat{1&0&0 \\ 0& \HF a_2 & b_2 \\0& -\overline{b_2} & \overline{a_2} }
\end{equation}
where $a_i,\, b_i \in \NC$ for $i=1,\,2$  satisfy $|a_i|^2 + |b_i|^2=1$ and $\zl\in \NC$ is of modulus $1$; $|\zl|=1$. 
Clearly 
\begin{equation}
g= \pmat{1&0&0 \\0 &\HF a_1 & b_1 \\0& -\overline{b_1} & \overline{a_1} }
\pmat{\zl^2\cos{\za}& -\sin{\za} &0 \\ \sin{\za}& \zl^{-2}\cos{\za}&0 \\ 0 &0 & 1}
\pmat{1&0&0 \\0 &\HF a_2 & b_2 \\0& -\overline{b_2} & \overline{a_2} }
\label{factor}
\end{equation}
\end{coro}

We finish this section by presenting the explicit orthogonal bases for the Lie algebras $\su{2}$  and $\su{3}$, which will be used in our future discussion of the $\C{P^1}$ and $\C{P^2}$ models. For uniformity we use the inner product 
\begin{equation}
\label{pr}
(X,Y) = -\frac{1}{2} \, \Tr(X  Y),
\end{equation}
where the Killing form is given by the formula 
\[
B(X,\,Y)= \begin{cases}
4\Tr(XY), & N=1\\
6\Tr(XY), & N=2 
\end{cases}
\qquad  X,\,Y\in \ali{su}{N+1}. 
\]
The orthogonal basis with respect to the Killing form is given for the case $N=1$ by the  matrices $-i\zs_j$, for $j=1,\,2,\,3$,  where $\zs_j$ denote the Pauli matrices
\begin{equation}\label{Pauli}
\zs_1=\pmat{0&1\\1&0},\quad \zs_2=\pmat{0&-i\\i&0},\quad \zs_3=\pmat{1&0\\0&-1}.
\end{equation}
Now, $\ali{su}{3}$ is eight-dimensional (over $\NR$) and consists of matrices of the form
\[
\pmat{ia_1 &   z_0 & z_1\\
-\overline{z}_0 &  ia_2 & z_2\\
-\overline{z}_1 &  -\overline{z}_2 & ia_3},
\qquad \text{with}\quad z_0,\,z_1,\,z_2\in \NC,\, a_1,\,a_2,\,a_3\in \NR,\,
a_1+a_2+a_3=0. 
\]
For this case we can choose a basis adapted to the decomposition $\mathfrak{g}=\mathfrak{k}_0 \oplus\mathfrak{p}$, where the isotropy Lie subalgebra $\mathfrak{k}_0$ is given by 
\begin{equation}
\mathfrak{k}_0 = \left\{\pmat{ia_0 & 0 & 0\\0& ia_1 & -\overline{z}\\0& z &ia_2}
\mid a_0+a_1+a_2=0,\, a_j\in \NR, z\in \NC\right\}. 
\end{equation}
and may be further decomposed as  
\begin{equation}\label{k-decomposition}
\mathfrak{k}_0 = \ali{u}{1}\oplus \ali{u}{2}.
\end{equation}
Accordingly, as its orthogonal basis, we take $\{S_j\mid j=1,\,\ldots,\,4\}$, where   
\begin{equation} \label{S_4}
S_j=\pmat{0 & 0\\ 0 & -i\zs_j},\qquad j=1,\,2,\,3,\qquad  
S_4= 
\pmat{-2i &0&0\\0&i&0\\0&0&i}.
\end{equation} 
Its orthogonal complement $\mathfrak{p}$ consists then of matrices defined in \Ref{eqn:Z(x)} which take the form
\[
Z(x)=\pmat{0& -\overline{x}_1&-\overline{x}_2\\ 
x_1 & 0 & 0 \\ x_2 &0 & 0},\qquad 
\text{where}\quad x = \pmat{x_1 \\ x_2} \in \NC^2 
\]
We suplement the above defined basis of $\mathfrak{k}_0$ by the following basis for $\mathfrak{p}$
\begin{gather}
S_5 = Z(e_1)= \pmat{0&-1&0\\1&0&0 \\0&0&0},\quad 
S_6= Z(e_2)= \pmat{0&0& -1\\0& 0&0 \\1& 0&0},\\
S_7 = Z(ie_1)= \pmat{0& i& 0 \\ i& 0&0 \\ 0 &0&0},\quad 
S_8=Z(ie_2)= \pmat{0&0& i\\0& 0&0 \\i&0&0}. \label{S_8}
\end{gather}
and observe that these bases are orthogonal with respect to the Killing form. Thus the matrices $\{S_1,\,\ldots,\,S_8\}$ form an orthogonal basis for $\ali{su}{3}$. 

\section{The structural equations of surfaces immersed in $\su{N+1}$ algebras obtained through $\cpn$ models} 
\label{sec:EQUATIONS}
In order to investigate immersions defined by means of solutions of the $\cpn$ models and, in particular, to envisage the moving frames and the corresponding Gauss-Weingarten and the
Gauss-Codazzi-Ricci equations, it is convenient to exploit the Euclidean structure of the $\su{N+1}$ Lie algebras leading to an identification $\mathbb{R}^{N(N+2)} \simeq \su{N+1}$ described in the previous section. 

Let us assume that the matrix $\mathbb{K}$ given by \Ref{1_14}
is constructed from a solution $f$ of the Euler--Lagrange equation
\Ref{1_11} defined on an open connected and simply connected set $\Omega \subset \C$.  The conservation law \Ref{1_16} then holds and implies that the matrix-valued 1-form
\begin{equation}
\label{3_2} 
\d X = i (\mathbb{K}^{\dagger} \, \d\xi + \mathbb{K}\, \d\bar{\xi})= 
i(\mathbb{K}^{\dagger} +\mathbb{K})\d\xi^1 -(\mathbb{K}^{\dagger} -\mathbb{K})\d\xi^2 
\end{equation}
is closed $(\d (\d X)=0)$ and takes values in the Lie algebra
$\su{N+1}$ of antihermitian matrices. By decomposing $\d X$ into the real and imaginary parts we write 
\begin{equation}
\label{3_4} \d X = \d X^1+i \d X^2,
\end{equation}
where the $1$-forms $\d X^1$ and $\d X^2$ with values in $\lie{sl}{N+1}{\NR}$ are anti-symmetric and symmetric, respectively, i.e. 
$$
(\d X^1)^T = -\d X^1 \, , \quad (\d X^2)^T = \d X^2,
$$
with the superscript $T$ denoting the transposition. From the closedness of the 1-form $\d X$ it follows  that the integral
\begin{equation}
\label{3_7} i \int_{\gamma} (\mathbb{K}^{\dagger} \d\xi +
\mathbb{K} \d\bar{\xi}) = X(\xi,\bar{\xi})
\end{equation}
locally depends only on the end points of the curve $\gamma$ (i.e. it is
locally independent of the trajectory in the complex plane $\C$).
The integral defines a mapping
\begin{equation}
\label{3_8} X: \Omega \ni (\xi,\bar{\xi}) \mapsto X(\xi,\bar{\xi}) \in \su{N+1}\simeq \mathbb{R}^{N(N+2)}
\end{equation}
which will be called the generalized Weierstrass formula for the immersion.  The
tangent vectors of this immersion, by virtue of \Ref{3_2}, are
\begin{equation}\label{eqn:real_der}
  \p_1 X= i(\mathbb{K}^{\dagger} +\mathbb{K}),\qquad 
\p_2 X= -(\mathbb{K}^{\dagger} -\mathbb{K})
\end{equation}
and the complex tangent vectors are
\begin{equation}
\label{3_10} \p X =i\mathbb{K}^{\dagger} \, , \quad 
\bp X = i \mathbb{K}.
\end{equation}

Hence a surface $\mathbb{F}$ associated with the $\cpn$ model \Ref{1_11} by means of the immersion \Ref{3_8}  satisfies the following 
\begin{theo}[Metric form]\label{theo:metric_form}\hfill\\
Components of the metric form induced on $\mathbb{F}$ by the Euclidean structure in $\su{N+1}$ defined by the negative of the Killing form \Ref{eqn:Killing} are given by 
\begin{align}
g_{11}  = (\p_1X,\, \p_1X)  = &\hphantom{{}-{}}2(N+1)\tr(\mathbb{K}^{\dagger} +\mathbb{K})^2 \nonumber \\
g_{22}  = (\p_2X,\, \p_2X ) = &-2(N+1)\tr(\mathbb{K}^{\dagger} -\mathbb{K})^2  \\
g_{12}  = (\p_1X,\, \p_2X ) = &\hphantom{{}-{}} 2i(N+1)\tr[(\mathbb{K}^{\dagger} +\mathbb{K})(\mathbb{K}^{\dagger} -\mathbb{K})]. \nonumber   
\end{align}
The components of the metric form with respect to the complex tangent vectors are given by the following expressions  
\begin{gather*}
\label{3_11}
g_{\xi,\xi}         = (\p X, \p X) = J, \qquad 
g_{\bar\xi,\bar\xi} =  (\bp X, \bp X) = \bar{J} ,\\
g_{\xi,\bar\xi} =  g_{\bar\xi,\xi} = (\p X, \bp X) = q,   
\end{gather*}
where $J$ and $\bar{J}$ are functions defined by \Ref{1_18}  
and $q$ is a non-negative (real valued) function given by  
\begin{equation*}
q  = \frac{1}{\fd f} \bp \fd P \p f \ge 0.
\end{equation*}
Consequently the first fundamental form $I$ of the surface $\mathbb{F}$ is
given with respect to the complex coordinates $\xi,\,\bar\xi$  by
\begin{equation}
\label{3_14} I= J\, \d\xi^2 + 2q \, \d\xi \,
\d\bar{\xi} + \bar{J}\,\d\bar{\xi}^2. 
\end{equation}
\end{theo}
In Section \ref{sec:CP1} we compute explicitely coefficients  of the metric form in the case of $\C P^1$ model. 

As usual, we denote 
$$
g= g_{\xi,\xi}g_{\bar\xi,\bar\xi}-g_{\xi,\bar\xi}^2= |J|^2-q^2
$$ 
the determinant of the metric form. It is known that the Gaussian curvature of the surface $\mathbb{F}$ with respect to the induced metric is given by 
\begin{equation}
K = \frac{1}{2\sqrt{g}}\bp\left[\frac{1}{\sqrt{g}}(-2\p q +q\p(\ln J))\right].
\label{K4_9}
\end{equation}
The quantity $J\, \d \xi^2$ defined on $\mathbb{F}$, called Hopf differential, is
invariant with respect to conformal changes of coordinates.  We use this
freedom to simplify the corresponding equations.

Our next task is to determine a moving frame on the surface $\mathbb{F}$ and to write the
corresponding Gauss-Weingarten equations expressed in terms of a solution $f$ satisfying the $\cpn$ sigma model equations \Ref{1_11}.  Using the Gramm-Schmidt orthogonalization procedure
to construct and write explicitly expressions for the normals
$\eta_k$ to a given surface in $\mathbb{R}^{N(N+2)}$ can, in general, be a
 complicated task. An alternative way we propose here involves the use of the isomorphism of 
$\mathbb{R}^{N(N+2)}$ with the Lie algebra $\su{N+1}$. In this representation,
the equations for a moving frame on the surface can be written in
terms of $(N+1)$ by $(N+1)$ skew-hermitian matrices. To simplify, in the following calculations we suppress the normalization factor $2(N+1)$ in the definition of the inner product in $\su{N+1}$ --- cf. \Ref{eqn:Killing}. We introduce real normal vectors $\eta_3, \ldots, \eta_{N(N+2)}$ to the surface  $\mathbb{F}$  and consider the frame
\begin{equation*}
\eta=\left(\eta_1=\p X, \eta_2 = \bp X, \eta_3, \ldots ,
\eta_{N(N+2)}\right)^T
\end{equation*}
with components satisfying the following conditions
\begin{equation}
\label{3_17}
\begin{split}
&(\p X, \eta_k)=0 \, , \quad (\bp X, \eta_k)=0, \\
&(\eta_j,\eta_k)=\delta_{jk} \, , \quad j,k=3,\ldots, N(N+2).
\end{split}
\end{equation}
Next we define 
\begin{equation}
J_k = \tr{(\partial^2X\cdot\eta_k)}\mbox{,}\quad H_k = \tr{(\partial\bar{\partial}X\cdot\eta_k)}\mbox{.}
\end{equation}
Now we can formulate the following
\begin{theo}[The structural equations]{}\hfill\label{prop3}\\
For any solution $f$ of the $\cpn$ sigma model equations \Ref{1_11}, such that the determinant of the induced metric $g$ is nonzero in some neighborhood of a regular point
$p=(\xi_0,\bar{\xi}_0)$ in $\C$, there exists in this neighborhood
a moving frame $\eta$ on this surface which satisfies the following
Gauss--Weingarten equations
\begin{equation}
\label{3_18}
\p \eta_i = A_{il}\eta_l \quad , \bp \eta_i = B_{il}
\eta_l, \quad i,l=1,\ldots, N(N+2)
\end{equation}
where the $N(N+2)$ by $N(N+2)$ matrices $A$ and $B$ have the form
\begin{equation}
\label{3_19}
A=\left(%
\begin{array}{ccccc}
  a_{1,1} & a_{1,2} & J_3 & \cdots & J_{N(N+2)} \\
  0 & 0 & H_3 & \cdots & H_{N(N+2)} \\
  \alpha_{1,3} & \beta_{1,3} & 0 & \cdots & S_{3,N(N+2)} \\
  \vdots & \vdots & \vdots &  & \vdots \\
  \alpha_{1,N(N+2)} & \beta_{1,N(N+2)} & -S_{3,N(N+2)} & \cdots & 0 \\
\end{array}%
\right), 
\end{equation}
and
\begin{equation}
\label{3_20}
B=\left(%
\begin{array}{ccccc}
  0 & 0 & H_3 & \cdots & H_{N(N+2)} \\
  a_{2,1} & a_{2,2} & \bar{J}_3 & \cdots & \bar{J}_{N(N+2)} \\
  \alpha_{2,3} & \beta_{2,3} & 0 &  & \bar{S}_{3,N(N+2)} \\
  \vdots & \vdots & \vdots &  & \vdots \\
  \alpha_{2,N(N+2)} & \beta_{2,N(N+2)} & -\bar{S}_{3,N(N+2)} & \cdots & 0 \\
\end{array}%
\right). 
\end{equation}
The elements of $A$ and $B$ take the form
\begin{eqnarray}
\label{3_21} & &S_{jk}+S_{kj}=0 \, , \quad
\bar{S}_{jk}+\bar{S}_{kj}=0 \,
, \quad j \neq k\nonumber\\
& & \alpha_{1,j} = \frac{1}{g} (H_j\, g_{\xi,\bar\xi} - J_j \,
g_{\bar\xi,\bar\xi}) \, , \quad \beta_{1,j}=\frac{1}{g} (J_j \,
g_{\xi,\bar\xi} - H_j \, g_{\xi,\xi}),\\
& & \alpha_{2,j} = \frac{1}{g} (\bar{J}\,
g_{\xi,\bar\xi}-H_j\, g_{\bar\xi,\bar\xi}) \, , \quad \beta_{2,j} = \frac{1}{g} (H_j \, g_{\xi,\bar\xi} - \bar{J}_j \, g_{\xi,\xi}),\notag
\end{eqnarray}
where
\begin{equation}
\label{3_22}
\begin{split}
& a_{1,1}  =  -a_{2,2}= \frac{1}{g} \mathrm{Re} \left\{\frac{1}{\fd f}
\left(\bar{J} \p \fd + g_{\xi,\bar\xi} \bp \fd \right) P \p^2 f \right.
\\ &  \left. \qquad\qquad -\frac{2\p \fd f}{(\fd f)^2}\left(\bp \fd P \p f\right) g_{\xi,\bar\xi} -
\frac{2 \fd \p f}{\fd f} |J|^2\right\}, \\
& a_{1,2}  =  \frac{1}{g} \mathrm{Re} \left\{\frac{-1}{\fd f}
\left(J \bp \fd + g_{\xi,\bar\xi} \p \fd \right) P \p^2 f \right. \\
& \left. \qquad\qquad +\frac{2 \p \fd f}{(\fd f)^2} (\bp \fd P \p
f) \bar{J} + \frac{2\fd \p f}{\fd f} J\, g_{\xi,\bar\xi}\right\},\\
& a_{2,1} = \frac{1}{g} \mathrm{Re} \left\{ \frac{1}{\fd f}
\left(J \bp \fd + g_{\xi,\bar\xi} \p \fd \right) P \bp^2 f \right. \\
& \left. \qquad\qquad -\frac{2\bp \fd f}{(\fd f)^2} \left(\p \fd P
\bp f\right) g_{\xi,\bar\xi} - \frac{2\fd \bp f}{\fd f} |J|^2\right\}.\\
\end{split}
\end{equation}
The Gauss--Codazzi--Ricci equations are given by
\begin{equation}
\label{3_23} \bp A - \p B + [A,B]=0
\end{equation}
and coincide with the equations of the $\cpn$ sigma model \Ref{1_11}.
\end{theo}
We note that the elements $a_{i,j}$ are the usual Christoffel symbols.

{\bf Proof.} Note that for any solution of the $\cpn$ equations
\Ref{1_11} the matrices $\p X$ and $\bp X$ are defined by
\Ref{3_10}.  As can be checked by a straightforward computation
using \Ref{1_11}, the mixed derivatives $\p \bp
X$ and $\bp \p X$ coincide and are normal to the surface
\begin{eqnarray}
\p\bp X & = & [\p P, \bp P]\nonumber\\
        & = & \frac{1}{\fd f} (P \p f \otimes \bp \fd P - P \bp f
\otimes \p \fd P) \nonumber\\
        &   &  \qquad + \frac{1}{(\fd f)^2} (\p \fd P \bp f - \bp
\fd P \p f) f \otimes \fd\nonumber\\
        & = & -[\bp P, \p P] = \bp \p X .
\end{eqnarray}
Combining this equation with the $\cpn$ equations \Ref{1_11},
expressed in terms of the projector P, we obtain
\begin{equation}
\begin{split}
&\tr{(\p\bp X \cdot \p X)} = \tr{([\p P, \bp P] \cdot [\p P,P])} =
0,\\
&\tr{(\p \bp X \cdot \bp X)} = \tr{([\p P, \bp P] \cdot [\bp P,
P])} = 0.
\end{split}
\end{equation}
As a direct consequence of differentiation of the normals
\Ref{3_17} we get
\begin{equation}
\begin{split}
& (\p \eta_j, \eta_k) + (\p \eta_k, \eta_j) = S_{jk}+S_{kj} = 0,\\
& (\bp \eta_j , \eta_k) + (\bp \eta_k, \eta_j) = \bar{S}_{jk} +
\bar{S}_{kj} = 0 ,\quad j \neq k
\end{split}
\end{equation}
and
\begin{equation}
\label{3_27}
\begin{split}
 (\bp \eta_j, \p X) + (\eta_j, \p \bp X) &= 0,\\
 (\bp \eta_j, \bp X) + (\eta_j, \bp^2 X) &= 0,\\
 (\p \eta_j, \p X) + (\eta_j, \p^2 X) &= 0,\\
 (\p \eta_j, \bp X) + (\eta_j, \bp \p X) &= 0.
\end{split}
\end{equation}
Using the expressions \Ref{3_18}-\Ref{3_20} and \Ref{3_11})
we come to the following set of linear equations
\begin{equation}
\begin{split}
&  g_{\xi,\bar\xi}\alpha_{1,j}+g_{\bar\xi,\bar\xi}\beta_{1,j} + J_j=0, \\
&  g_{\xi,\xi}\alpha_{1,j} + g_{\xi,\bar\xi}\beta_{1,j} + H_j = 0, \quad j=3,
\ldots, N(N+2)\\
&  g_{\xi,\bar\xi} \alpha_{2,j} + g_{\bar\xi,\bar\xi}\beta_{2,j} + \bar{J}_j = 0,\\
&  g_{\xi,\xi}\alpha_{2,j} + g_{\xi,\bar\xi}\beta_{2,j} + H_j = 0,
\end{split}
\end{equation}
which allow us to determine elements $\alpha_{i,j}$ and $\beta_{i,j}$ in terms of the coefficients of the metric form, $H_j$ and $J_j$.  As can be
easily calculated, they take the form \Ref{3_21} postulated in
Proposition \ref{prop3}.  The second derivatives $\p^2 X$ and
$\bp^2 X$ are
\begin{equation}
\label{3_30}
\begin{split}
\p^2 X  &=  
\frac{1}{\fd f} (P \p^2 f \otimes \fd - f
\otimes \p^2 f P)\\
 &+ \frac{2}{(\fd f)^2} \left[(\p \fd \cdot f) f
\otimes
\p \fd P - (\fd \cdot \p f) P \p f \otimes \fd \right],\\
\bp^2 X  &=  
\frac{1}{\fd f} (f \otimes \bp^2 \fd P - P
\bp^2 f \otimes \fd)\\
 &+ \frac{2}{(\fd f)^2} \left[(\fd \cdot \bp f)
P \bp f \otimes \fd - (\bp \fd \cdot f) f \otimes \bp \fd
P\right].
\end{split}
\end{equation}
Let us observe that the following traces (and their complex
conjugates) vanish:
\begin{equation}
\begin{split}
& \tr{\left(\left(\p^2 X -a_{1,1} \p X -a_{1,2}\bp X \right) \cdot
\p X \right)} = 0,\\
& \tr{\left(\left(\p^2 X - a_{2,1} \p X -a_{2,2} \bp X \right) \cdot
\bp X\right)} = 0.
\end{split}
\end{equation}
This means that the vectors corresponding to the matrices $(\p^2 X -
a_{i,1} \p X -a_{i,2} \bp X)$ and $(\bp^2 X - a_{i,1}\p X -a_{i,2} \bp
X)\, , i=1,2$ are normal to the surface determined by
\Ref{3_18}. Substituting \Ref{3_30} into equations
\Ref{3_27} and solving the obtained linear systems we can
determine  $a_{i,l}$, $i,l=1,2$ which
prove to be of the form \Ref{3_22}.\\
\indent Finally, the Gauss--Codazzi--Ricci (GCR) equations are the
necessary and sufficient conditions for a local existence of a
surface and are the compatibility conditions of the
Gauss-Weingarten equations.  In our case the GCR equations
coincide with the $\cpn$ sigma model equations \Ref{1_11} and are
given in a matrix form by \Ref{3_23}.  So, with any solution $f$
of the $\cpn$ model we can associate a surface defined by 
\Ref{3_7}.  {\bf Q.E.D.}\\

Making use of the expressions for the second derivatives of $X$,  
the induced metric and for the elements $a_{i,l}$ appearing in
matrices $A$ and $B$, we can write explicitly the second
fundamental form of the surface in terms of the model 
\begin{equation}
\label{3_31}
\begin{split}
\mathbf{II} & = (\p^2 X)^{\bot} \d \xi^2 + 2(\p \bp X)^{\bot} \d
\xi \d \bar{\xi} + (\bp^2 X)^{\bot} \d \bar{\xi}^2\\
& = (\p^2 X -a_{1,1}\p X -a_{1,2} \bp X) \d \xi^2 + 2(\p \bp X) \d
\xi \d \bar{\xi} \\
& \quad + (\bp^2 X -a_{2,1}\p X -a_{2,2} \bp X)\d \bar{\xi}^2,
\end{split}
\end{equation}
where the symbol $\bot$ denotes the normal part of matrices
$\p_i\p_jX$ and the indices $i,j$ stand for $\xi$ or $\bar{\xi}$.
The quantities $a_{i,l}$ are given by \Ref{3_22}. \\
The mean curvature vector can also be expressed in terms of the model as
follows:
\begin{equation}
\label{3_32}
\begin{split}
\mathbf{H} & =\frac{1}{g} (g_{\xi,\xi} (\bp^2 X)^{\bot} - 2g_{\xi,\bar\xi}
(\p\bp X)^{\bot} +g_{\bar\xi,\bar\xi} (\p^2 X)^{\bot})\\
& = \frac{1}{g} \left\{ g_{\bar\xi,\bar\xi}[\p^2X-a_{1,1} \p X - a_{1,2}
\bp X] -2g_{\xi,\bar\xi} \p \bp X \right.\\
& \left. \qquad + g_{\xi,\xi}[\bp^2 X -a_{2,1} \p X -a_{2,2} \bp
X]\right\}.
\end{split}
\end{equation}
The Willmore functional of a surface has the form
\begin{equation}
W=\int_{\Omega} |\mathbf{H}|^2 \sqrt{g} \,\d \xi\, \d \bar{\xi}.
\end{equation}
When a solution $f$ satisfying the $CP^N$ model \Ref{1_11} is defined on the whole Riemannian sphere ${S}^2$ then the integral
\begin{eqnarray}
\label{3_35} Q & = & \frac{i}{8\pi} \int_{S^2} \tr{(P \cdot [\p P, \bp P])}
\d \xi \d \bar{\xi} \nonumber\\
  & = & \frac{1}{8\pi} \int_{S^2} \frac{1}{\fd f} (\p \fd P \bp
f - \bp \fd P \p f) \d \xi \d \bar{\xi}
\end{eqnarray}
is an invariant of the surface and is known as the topological charge of
the model.  If the integral \Ref{3_35} exists then it is an 
integer which characterizes globally the surface under
consideration.\\
Summarizing, we can now state the following analog of the Bonnet theorem, cf. \cite{Wil}.
\begin{coro} For the complex-valued function $f$ satisfying
the $\cpn$ sigma model equations \Ref{1_11}, the generalized Weierstrass formula for immersion \Ref{3_8}, i.e.
\begin{equation}
X:\Omega\ni (\xi,\bar{\xi}) \rightarrow X(\xi,\bar{\xi}) = i\int_{\gamma}[\partial P,P]d\xi + [\bar{\partial}P,P]d\bar{\xi}\mbox{,}
\end{equation}
 describes a surface in $\su{N+1}$. This surface is determined 
by its first and second fundamental forms \Ref{3_14} and \Ref{3_31} uniquely up to Euclidean motions. 
\end{coro}

Finally, it is worth noting that the method described above may be of use in the study of the elliptic periodic two-dimensional Toda lattice ($2DTL$) which is related to surfaces immersed in $\su{N+1}$ Lie algebra \cite{GU2}{}. The equations of $2DTL$ can be written in a matrix form as the zero curvature equations $\bp A - \p B = [A,\,B]$, formally identical with the  Gauss--Codazzi--Ricci equation \Ref{3_23}, where the two $(N+1)$ by $(N+1)$ matrices  $A$ and $B$ are defined as follows
\begin{equation}
\label{4.30}
A=-B^\dagger =\pmat{\p u_0 & 0 &0 &\ldots &0& U_{0,N}\\
U_{1,0}& \p u_1  & 0 & \ldots & 0& 0\\
0 &U_{2,1} & \p u_2  & \ldots &0&0\\
\vdots &\vdots& \vdots & & \vdots &\vdots\\
0 & 0 &0 &\ldots & \p u_{N-1}& 0\\
0 & 0 &0 &\ldots &  U_{N,N-1}& \p u_N
},
\end{equation}
where for $i,j=0,\ldots,N$, we set $u_i:\C\to \NR$,   
$u_0+\ldots +u_N=0$ and $U_{i,j}=\exp(u_i-u_j)$. It is known \cite{GU2} that the 
zero-curvature equation \Ref{3_23} for matrices \Ref{4.30} implies the existence of a complex valued function $F:\C\to \gr{SU}{N+1}$ such that 
\begin{equation}
F^{-1}\p F = A,\qquad F^{-1}\bp F= B
\label{4.31}
\end{equation}
So, according to the Proposition \ref{prop3}, we can identify \Ref{4.31} with the complex tangent vectors \Ref{3_10} of the immersion \Ref{3_8}. Hence the $2DTL$ equations can be viewed as being associated with the specific form \Ref{1_13} of the $\cpn$ model. Establishing this link could be useful for determining certain geometric characteristics of surfaces corresponding to the elliptic periodic two-dimensional Toda lattice, but this point will not be considered here. 

\section{Immersions into the Lie algebra $\su{2}$ arising from the $\C P^1$ model}
\label{sec:CP1}

In this section we sketch an application of the techniques developed in the previous sections to the case of $\C P^1$ sigma model. This allows us to put the results obtained in the earlier publications\cite{GZ01,GZ,GRZA} in a broader perspective, as well as to point out some further geometrical properties of surfaces obtained from this model. 

The fields of the $\C P^1$ model in the notation of the Section \ref{sec:SUN} are given by $[z]=[f_0,\,f_1]$, where $z=(f_0,\,f_1)\in S^3\subset \C^2$, but it is customary to replace here the homogeneous coordinates $(f_0,\,f_1)$ by the affine coordinate $W=f_1/f_0$. The Euler--Lagrange equation \Ref{1_11} reduces to   
\begin{equation}
\label{1_22'}
\partial \bar{\partial} W - \frac{2\bw}{1+|W|^2} \, \partial W
\,\bp W=0
\end{equation}
and the matrix $\mathbb{K}$ of the eqn. \Ref{1_14} is then given by
\begin{equation}
\mathbb{K} = \frac{1}{(1+|W|^2)^2} \left( \begin{array}{cc}
\bw \bp W - W \bp \bw & \bp \bw + \bw^2 \bp W \\
-\bp W - W^2\bp \bw & W \bp \bw - \bw\bp W \end{array} \right).
\end{equation}
 
The $1$-form $\d X=i(\mathbb{K}^\dagger\d \xi+\mathbb{K}\d \bar{\xi})$ of the eqn. \Ref{3_2} defines an immersion into $\su{2}$ --- to pass to the Euclidean space $\NR^3$ we first compute its real and imaginary parts and obtain  
\begin{eqnarray*}
\textstyle \d X^1 &=& \textstyle \frac{i}{2(1+|W|^2)^2}
\left[\left(
\begin{array}{cc} 
 0 & 
-\p \bw-\bw^2\p W - \p W - W^2 \p
\bw \\ 
\p W+W^2\p \bw + \p \bw + \bw^2 \p W &
 0
\end{array} \right)\d\xi \right.\\
& & \left.+ \left( \begin{array}{cc} 
 0 & 
 \bp \bw +\bw^2
\bp W + \bp W +W^2 \bp \bw \\ 
 -\bp W -W^2\bp \bw-\bp \bw-\bw^2
\bp W & 
 0 \end{array} \right) \d \bar{\xi} \right],
\end{eqnarray*}
and
\begin{eqnarray*}
\textstyle \d X^2 &=& \textstyle \frac{1}{2(1+|W|^2)^2}
\left[\left(
\begin{array}{cc} 
 2(W\p \bw - \bw\p W) & 
 -\p \bw-\bw^2\p W + \p W + W^2 \p
\bw \\ 
 \p W+W^2\p \bw - \p \bw - \bw^2 \p W & 
 -2(W \p \bw
- \bw\p W)
\end{array} \right)\d\xi \right.\\
& & \left.+ \left( \begin{array}{cc} 
 2(\bw\bp W - W \bp
\bw) & 
 \bp \bw +\bw^2 \bp W - \bp W - W^2 \bp \bw \\
 -\bp W -W^2\bp \bw+\bp \bw+\bw^2 \bp W & 
 -2(\bw\bp W -
W\bp \bw)
\end{array} \right) \d \bar{\xi} \right].
\end{eqnarray*}
These in turn can be expressed in terms of the Pauli matrices $\sigma_i$ (cf. eqn. \Ref{Pauli}) as follows   
\begin{equation*}
\d X^1 = i\d X_1 \sigma_2 \, , \quad \d X^2 = \d X_2\sigma_1 + \d
X_3 \sigma_3, 
\end{equation*}
where the real-valued $1$-forms $\d X_i,\ i=1,2,3$, are given by  
\begin{equation}
\label{4_3}
\begin{split}
\d X_1 & = 
\frac{i}{2(1+|W|^2)^2}
\Big{\{}-\left[\p \bw + \bw^2\p W + \p W +W^2\p \bw\right] \d \xi\\
& +
\left[\bp \bw +\bw^2\bp W + \bp W +W^2 \bp \bw\right] \d \bar{\xi}
\Big{\}},\\
\d X_2 & =  
\frac{1}{2(1+|W|^2)^2} \Big{\{}\left[-\p \bw -
\bw^2\p W + \p W +W^2\p \bw\right] \d \xi \\
& + \left[\bp \bw
+\bw^2\bp W - \bp W -W^2 \bp \bw\right] \d \bar{\xi}
\Big{\}},\\
\d X_3 & =  
\frac{1}{(1+|W|^2)^2} \left\{\left[W \p\bw -
\bw \p W\right] \d \xi + \left[\bar{W}\bp W - W \bp \bw\right]\d
\bar{\xi}\right\},
\end{split}
\end{equation}
which is the generalized Weierstrass formula for an immersion into $\NR^3\simeq \su{2}$. An interested reader may check that the equations \Ref{4_3} yield the classical Weierstrass formula \Ref{minimal} under the substitution $W=f_1/f_0$. 

Starting from a particular solution $W$ of the $\C{P}^1$ sigma model equation \Ref{1_22'} one constructs an immersion in $\mathbb{R}^3$ by the use of the formulae \Ref{4_3}.  The following is now readily obtained form the Proposition \ref{theo:metric_form}. 
\begin{coro} 
For the immersion given by equations \Ref{4_3} the coefficients of the induced metric are given by the following expressions 
\begin{align}\label{real_metric} 
g_{11}  =  & \frac{|\p W|^2+|\bp W|^2+ |\p W - \bp W|^2 }{(1+|W|^2)^2},\nonumber \\
g_{22}  =  & \frac{|\p W|^2+|\bp W|^2+ |\p W + \bp W|^2 }{(1+|W|^2)^2}, \\
g_{12}  =  & \frac{2\Im(\p W \p \bw)}{(1+|W|^2)^2}.  \nonumber		
\end{align}
The complex form of the induced metric in this case is given by 
\begin{equation}\label{metric}
g_{\xi,\xi}         = -\frac{ \p W \p \bw}{(1+|W|^2)^2}, \quad    
g_{\bar\xi,\bar\xi} =  -\frac{\bp W \bp \bw}{(1+|W|^2)^2}, \quad 
g_{\bar\xi,\xi}     =  \hphantom{-} \frac{|\p W|^2+ |\bp W|^2}{(1+|W|^2)^2}.   
\end{equation}
\end{coro}
For solutions of \Ref{1_22'} which are defined over $S^2$, the function
$W$ can be only holomorphic or antiholomorphic (cf. \cite{WZ}). For holomorphic $W$   
equations \Ref{metric} reduce to 
\begin{equation}
g_{\xi,\xi}=g_{\bar\xi,\bar\xi}=J=0 \, , \qquad g_{\bar\xi,\xi} = \frac{|\p W|^2}{(1+|W|^2)^2} 
\end{equation}
implying that the immersion is conformal. These relations, as shown earlier \cite{GRZA}{}, imply also $g_{\bar\xi,\bar\xi}=|Dz|^2$. 
The Gaussian curvature is $K= 1$, and the first and the second fundamental forms for the immersion are equal,  
\begin{equation}
II = I = {|\partial W|^2\over \left(1+|W|^2\right)^2}d\xi d\bar{\xi}.
\end{equation}
Moreover, as shown by K. Kenmotsu in \cite{Ken}{}, the function $W$ represents the complex Gauss map of the surface.  Geometrically, all that means that solutions of the $\mathbb{C}P^1$ model \Ref{1_22'} defined over $S^2$ parametrize the standardly immersed sphere in $\mathbb{R}^3$. This had already been shown in the case of instanton solutions of the $\gr{SO}{3}$ sigma model in \cite{FG}{}. 

\section{Surfaces immersed in the $\su{3}$ algebra}
\label{sec:su3}

Here we apply our considerations to the $\mathbb{C}P^2$ model for which we construct the associated immersion of a surface $\mathbb{F}$ in $\mathbb{R}^8$ and compute some of its geometric characteristics. For the case of $N=2$ we can pass from the representation 
$[z]=[f_0,f_1,f_2]$ to the inhomogeneous (affine) coordinates $W_1=f_1/f_0$ and $W_2=f_2/f_0$.     Now the Euler--Lagrange equations \Ref{1_11} take the form  
\begin{equation}
\label{1_23}
\begin{split}
& \partial \bar{\partial} W_1 - \frac{2\bar{W_1}}{A} \p W_1 \bp
W_1 - \frac{\bar{W_2}}{A}(\p W_1 \bp W_2+ \bp W_1 \p W_2) =0 , \\
& \p \bp W_2 - \frac{2\bar{W_2}}{A} \p W_2 \bp W_2 -
\frac{\bar{W_1}}{A} (\p W_1 \bp W_2 + \bp W_1 \p W_2) = 0,
\end{split}
\end{equation}
where
\begin{equation} 
A=1+|W_1|^2+|W_2|^2.\label{1_25'}
\end{equation}  

As was noted in \cite{GZ01}{}, the metric induced by the immersion $X(\xi,\bar{\xi})$ 
in $\mathbb{R}^8\simeq\su{3}$ is conformal for holomorphic solutions of the $\C P^2$  model  defined over $S^2$ and is then given by  
\begin{equation}
\label{4_2} g_{\xi,\xi}=g_{\bar\xi,\bar\xi}=0 , \quad  g_{\bar\xi,\xi} = 
\frac{|\p W_1|^2 + |\p W_2|^2 + |W_1 \p W_2 - W_2 \p W_1|^2}{A^2} .
\end{equation}
We shall write 
\begin{equation}\label{4_4}
  g_{\bar\xi,\xi} = e^{\frac{1}{2} (u+\bar{u})},
\end{equation}
where $u$ is a complex-valued function of $\xi,\,\bar{\xi} \in \C$ given by 
\begin{equation}
u+\bar{u} = \ln{\left\{\frac{1}{A^2} [ |\p W_1|^2 + |\p W_2|^2 +
|W_1\p W_2 - W_2\p W_1|^2]\right\}}.
\end{equation}

Under the above circumstances the following holds:
\begin{theo}[Structural equations for holomorphic $\C{P}^2$ model]{}\hfill\\
Any set of holomorphic solutions $(W_i,\bw_i)$, $i=1,2$ of the $\C P^2$ sigma
model equations \Ref{1_23} defined over $S^2$ such that the induced metric is nonzero in some neighbourhood of a regular point $p = (\xi_0,\bar{\xi}_0) \in \C$, 
determines a conformal parametrization of a surface $\mathbb{F}$ immersed in the $\su{3}$ Lie algebra. Its moving frame on $\mathbb{F}$ can be written in terms of $3$ by $3$ skew-hermitian matrices and is of the form
\begin{equation}
\label{4_5} \eta=(\p X, \bp X, \eta_1, \ldots, \eta_6)^{T},
\end{equation}
where the complex tangent vectors to the surface $\mathbb{F}$ are given by the following traceless matrices 
\begin{equation}
\label{4_6}
\begin{split} & 
 \p X = 
 -\frac{i}{A^2} \left(
\begin{array}{ccc} 
 -(\bw_1\p W_1+\bw_2 \p W_2) & 
 -\bw_1 (\bw_1\p W_1
+ \bw_2 \p W_2) & 
 -\bw_2(\bw_1\p W_1 + \bw_2 \p W_2) \\
 d_1 & 
 \bw_1d_1 & 
 \bw_2d_1 \\
 d_2 & 
\bw_1d_2 & 
\bw_2d_2
\end{array}
\right), \\ 
\noalign{\bigskip}
 & 
 \bp X = -\frac{i}{A^2} \left(
\begin{array}{ccc}
 W_1\bp \bw_1+W_2\bp \bw_2 & 
 -\bar{d_1} & 
 -\bar{d_2}\\
 W_1(W_1\bp \bw_1 + W_2\bp \bw_2) & 
 -W_1\bar{d_1} & 
 -W_1\bar{d_2} \\ 
 W_2(W_1\bp \bw_1 +W_2 \bp \bw_2) & 
-W_2\bar{d_1} & 
-W_2\bar{d_2} \end{array} \right), 
\end{split}
\end{equation}
and where we have defined 
\begin{equation}
\begin{split}
d_1 = (1+|W_2|^2)\partial W_1 - W_1\bar{W_2}\partial W_2,\\
d_2 = (1+|W_1|^2)\partial W_2 - \bar{W_1}W_2\partial W_1.
\end{split}
\label{dd}
\end{equation}
\end{theo}
\noindent{\bf Remark 2.} The explicit expressions for the complex normals to this surface immersed in $\su{3}$, in terms of $W_1$ and $W_2$, can be found in the Appendix to \cite{GSZ}. 

{\bf Proof.} Due to the normalization of the function $X$ (given
by equations \Ref{4_4} and \Ref{3_10}) we can express the
moving frame $\eta = (\p X, \bp X, \eta_1, \ldots , \eta_6)^{T}$
on a surface $\mathbb{F}$ in terms of the adjoint $\su{3}$
representation
\begin{equation}
\label{4_7}
\begin{cases}
\p X = e^{u/2} \zF^{-1} Y_{-} \zF, \\
\bp X = e^{\bar{u}/2} \zF^{-1} Y_{+} \zF,\\
\eta_i = \zF^{-1} S_{i+2} \zF , \quad i=1,\ldots,6
\end{cases}
\end{equation}
where
\begin{eqnarray}
Y_{-} = \frac{i}{2} (S_1 - i\,S_2) = \left(\begin{array}{ccc} 0 &
0 & 0 \\ 1 & 0 & 0 \\ 0 & 0 & 0 \end{array} \right) & , &  \quad
Y_{+} = \frac{i}{2} (S_1 + i\, S_2) = \left(\begin{array}{ccc} 0 &
1 & 0
\\ 0 & 0 & 0 \\ 0 & 0 & 0 \end{array} \right),\nonumber\\& &\\
(\zF^{-1} Y_{-} \zF)^{\dagger} & = & \zF^{-1} Y_{+} \zF.
\nonumber
\end{eqnarray}
Note that $\{Y_{-},Y_{+}\}$ span over $\mathbb{R}$, the same space
as $\{S_1,S_2\}$.  Using the polar decomposition of the $\SU{3}$
group given in Section \ref{sec:SUN}, cf. \Ref{eqn:Rowe_decomp}, a general $\SU{3}$ matrix $\zF$ can be decomposed into a product of three $\SU{2}$ factors. Performing
the multiplication in the expression \Ref{factor} and setting $\lambda=e^{i\phi/2}$ and $\alpha=t$, we obtain
\begin{equation}
\label{4_9} \zF = \left(\begin{array}{ccc} e^{i\varphi}\cos{t} &
-a_2\sin{t} & -b_2\sin{t} \\ a_1\sin{t} & a_1a_2
e^{-i\varphi}\cos{t}
-b_1\bar{b_2} & b_2a_1 e^{-i \varphi}\cos{t} +\bar{a_2}b_1 \\
-\bar{b_1}\sin{t} & -a_2\bar{b_1} e^{-i\varphi} \cos{t}
-\bar{a_1}\bar{b_2} & -\bar{b_1}b_2 e^{-i\varphi}
\cos{t}+\bar{a_1}\bar{a_2} \end{array} \right) \in \SU{3},
\end{equation}
where the complex-valued functions $a_i,\,b_i$ of $\xi,\,\bar{\xi}$
satisfy and $\varphi,\,t$ are real-valued functions of $\xi,\,\bar{\xi} \in
\C$.  The requirement that the parametrization of a surface $\mathbb{F}$ is conformal implies the following relations 
\begin{equation}
\begin{split}
& (\p X,\p X) = e^u \tr{(Y_{-})^2} = 0 ,\qquad  (\bp X, \bp X) =
e^{\bar{u}} \tr{(Y_{+})^2} = 0 , \\
& \quad \qquad (\p X, \bp X) = e^{\frac{1}{2} (u+\bar{u})}
\tr{(Y_{-} \cdot Y_{+} )} = e^{\frac{1}{2} (u+\bar{u})},
\end{split}
\end{equation}
and
\begin{equation}
\begin{split}
(\p X, \eta_i) = e^{u/2} \tr{(Y_{-} \cdot S_{i+2})} = 0,\\
(\bp X, \eta_i) = e^{\bar{u}/2} \tr{(Y_{+} \cdot S_{i+2})} = 0,\\
(\eta_j,\eta_k) = \tr{(S_{j+2} \cdot S_{k+2})} = \delta_{jk}.
\end{split}
\end{equation}

Now we have to determine the form of a 8-parameter representation of
the matrix $\zF$ in terms of $W_i, \bw_i$, compatible with the
$\C P^2$ sigma model \Ref{1_23}.  Using the $3$ by $3$ projector
matrix
\begin{equation}
\label{4_14} P = \ID_3- \frac{1}{A} \left( \begin{array}{ccc} 1 & W_1 & W_2 \\
\bw_1 & W_1 \bw_1 & \bw_1 W_2 \\ \bw_2 & W_1 \bw_2 & \bw_2 W_2
\end{array} \right), 
\end{equation}
we can write the Euler--Lagrange equations \Ref{1_23} in the form
of a conservation law \Ref{1_16} for the matrix
\begin{equation}
\begin{split}
\mathbb{K}=&\frac{1}{A} \left( \begin{array}{ccc} 0 & -\bp \bw_1 &
-\bp \bw_2 \\ \bp W_1 & \bw_1\bp W_1 - W_1 \bp \bw_1 & \bw_2 \bp
W_1 - W_1 \bp \bw_2 \\ \bp W_2 & \bw_1 \bp W_2 - W_2 \bp \bw_1 &
\bw_2 \bp W_2 - W_2 \bp \bw_2
\end{array} \right) \\&+ \frac{\bar{\rho}}{A^2} \left( \begin{array}{ccc} 1
& \bw_1 & \bw_2 \\ W_1 & |W_1|^2 & W_1\bw_2 \\ W_2 & \bw_1 W_2 &
|W_2|^2 \end{array} \right),
\label{KK}
\end{split}
\end{equation}
where we have defined the following expression
\begin{equation}
\rho = \bw_1 \p W_1 - W_1 \p \bw_1 + \bw_2 \p W_2 - W_2 \p \bw_2
\end{equation}
and the quantity $A$ is given by equation \Ref{1_25'}.
According to \Ref{1_14} and \Ref{3_10}, the matrices $\partial X$ and $\bar{\partial}X$ take the required form \Ref{4_6}.

Let us note that to satisfy the compatibility condition for
\Ref{4_7} (i.e. $\bp \p X = \p \bp X)$ it is sufficient, in view
of the conservation laws \Ref{1_13}, to postulate that the condition
\Ref{3_10} holds for the matrix $\mathbb{K}$ given by
\Ref{KK}.  So, we can determine functions $a_i,b_i \in \C$ and
$t,\varphi \in \mathbb{R}$, appearing in the matrix $\zF \in
\SU{3}$, in terms of $W_i$ and $\bw_i$,\, $i=1,2$.  
By a straightforward algebraic computation we get
\begin{equation}
\zF = \left(\begin{array}{ccc}  {e^{i \varphi}}{A^{-1/2}} &
\bw_1
{e^{i \varphi}}{A^{-1/2}} & \bw_2 {e^{i \varphi}}{A^{-1/2}} \\
\noalign{\medskip}
{ir^{-1} e^{i\varphi} (W_1 \p W_1 +W_2  \p W_2)}{A^{-1/2}} &
{-ir^{-1} e^{i \varphi} }{A^{-1/2}}\bar{d_1} & {-ir^{-1} e^{i \varphi} }{A^{-1/2}}\bar{d_2} \\
\noalign{\medskip}
{ir^{-1} (W_1\p W_2 - W_2 \p W_1)e^{-2i \varphi}} & {-ir^{-1} \p
W_2 e^{-2i \varphi}} & {ir^{-1} \p W_1e^{-2i \varphi}}
\end{array} \right),
\end{equation}
where we have used the notation introduced in \Ref{dd} and have set $r^2=A^2g_{\bar\xi,\xi}$.

Given the above form of the matrix $\zF$, the matrices $Y_{-}$,
$Y_{+}$ and the $S_{i+2}$,\, $i=1,\ldots, 6$ the moving frame
\Ref{4_7} adopts the required forms \Ref{4_5}-\Ref{4_6}. One
can check directly that it satisfies the Gauss-Weingarten
equations \Ref{3_18}.  In our case the corresponding GCR
equations, which are the compatibility conditions for
\Ref{3_18}, coincide with the $\C P^2$ sigma model equations
\Ref{1_23}.  Thus we have proved that any holomorphic solution of the $\C P^2$ 
model defined over $S^2$ gives a surface  conformally immersed in $\mathbb{R}^8$ with the moving frame given by \Ref{4_5} and \Ref{4_6}. {\bf Q.E.D}\\

\section{Examples and applications for the $\mathbb{C}P^2$ model}
\label{sec:CP2}

Based on the results of the previous sections we can
now construct certain classes of two-dimensional surfaces immersed
in $\mathbb{R}^8$. For this purpose we use the $\C P^2$ sigma
model defined over $S^2$.  For this model all solutions of the
Euler--Lagrange equations \Ref{1_23} are known \cite{WZ}{}. If we 
require the finiteness of the action\Ref{1_7} they
split into three classes : holomorphic (i.e. $W_i=W_i(\xi)$),
antiholomorphic (i.e. $W_i = W_i(\bar{\xi})$) and the mixed ones.
The latter ones can be determined from either the holomorphic or
antiholmorphic functions by the following procedure \cite{WZ}{}:\\
\indent Consider three arbitrary holomorphic functions
$g_i=g_i(\xi)$ and define for any pair of them the Wronskian functions
\begin{equation}\label{6_1}
G_{ij} = g_i \p g_j - g_j\p g_i \,,\quad  \bp g_i = 0\, , \quad
i=1,2,3
\end{equation}
Then one can check that the map $f=(f_1,\,f_2,\,f_3)$, where 
\begin{equation}
f_i = \sum_{k\neq i}^3 \bar{g}_k G_{ki} \,,\quad i=1,2,3
\end{equation}
is a solution of the $\C P^2$ sigma model, so called the mixed one, and hence  
the ratios 
\begin{equation}
W_1 = \frac{f_1}{f_3} \,,\quad W_2 = \frac{f_2}{f_3}
\label{5_3}
\end{equation}
satisfy equations \Ref{1_23}.

An alternative approach starts with any antiholomorphic functions
$\bar{g}_i = \bar{g}_i(\bar{\xi})$ and constructs functions $f_i$
and consequently $W_i$ as above but using $\bp$ instead of $\p$ in
the definition of $G_{ij}$.  It yields results which are
complementary to the ones obtained by the first approach.  Let us
note here that the requirement of finite action \Ref{1_7}
excludes solutions which admit Painlev{\'e} transcendent (i.e. all critical points are fixed independent of initial data), branch points or essential singularities.\\

\indent Now, let us discuss some classes of surfaces immersed in
$\su{3}$ algebra which can be obtained directly by applying the
Weierstrass representation \Ref{3_8}. For the $\C{P}^2$ model
\Ref{1_23}, the matrix $\mathbb{K}$ in terms of $W_i$ and
$\bar{W_i}$, $i=1,2$ is given by \Ref{KK}.
From the equation \Ref{3_2} we obtain for the real
and imaginary parts of the $1$-form $\mathrm{d}X$ the expressions 
\begin{equation}
\begin{split}
\mathrm{d}X^1 &=\frac{i}{2}
\left[(\mathbb{K}^{\dagger}-\bar{\mathbb{K}})\mathrm{d}\xi + (\mathbb{K} - \mathbb{K}^T) \mathrm{d}\bar{\xi}\,\right],\\
\mathrm{d}X^2 &=\frac{1}{2} \left[(\mathbb{K}^{\dagger} +
\bar{\mathbb{K}})\mathrm{d}\xi + (\mathbb{K}+\mathbb{K}^T)
\mathrm{d} \bar{\xi}\,\right].
\end{split}
\end{equation}
Clearly, the matrices $\d X^1$ and $\d X^2$ are antisymmetric and
symmetric, respectively, and hence can be decomposed in terms of the chosen basis in $\su{3}$ given by \Ref{S_4}--\Ref{S_8} as follows  
\begin{equation}
\label{5_7}
\begin{split}
\d X^1 &= \d X_2 S_2 + \d X_5 S_5 + \d X_6 S_6  ,\\
\d X^2 &= i\, \left[\d X_1 S_1 + \d X_3 S_3 + \d X_4 S_4 + \d X_7 S_7  + \d X_8 S_8 \right].
\end{split}
\end{equation}
As a result of the decomposition \Ref{5_7}, there exists eight real-valued functions $X_i (\xi, \bar{\xi})$ $i=1, \ldots, 8$ which determine the generalized Weierstrass representation of surfaces associated with the $\C P^2$ model \Ref{1_23}.  Considering the off-diagonal
entries of the matrices $\d X^1$ and $\d X^2$ we get
\begin{equation}
\label{5_8}
\begin{split}
\d X_1 &= [\frac{1}{A} (\p \bw_1 - \p W_1) +
\frac{\rho}{A^2}(\bw_1+W_1)]\d \xi + [\frac{1}{A}(\bp W_1 - \bp
\bw_1) + \frac{\bar{\rho}}{A^2} (\bw_1+W_1)]\d \bar{\xi},\\
\d X_2 &= -i\Big\{[\frac{1}{A} (\p W_1 + \p \bw_1) +
\frac{\rho}{A^2}(\bw_1 - W_1)]\d \xi \\
&  \quad    +[-\frac{1}{A} (\bp \bw_1 + \bp W_1) +
\frac{\bar{\rho}}{A^2} (\bw_1 - W_1) ] \d \bar{\xi}\Big\},\\
\d X_5 &= -i\Big\{[\frac{1}{A} (\p W_2 + \p \bw_2) +
\frac{\rho}{A^2}(\bw_2-W_2)]\d \xi \\
& \quad + [-\frac{1}{A}(\bp \bw_2 +\bp W_2) +
\frac{\bar{\rho}}{A^2} (\bw_2 - W_2)]\d \bar{\xi}\Big\},
\end{split}
\end{equation}
\begin{equation}
\begin{split}
\nonumber
\d X_6 &= -i\Big\{[\frac{1}{A} (W_1 \p \bw_2 - \bw_2 \p W_1 - W_2
\p \bw_1 + \bw_1 \p W_2) + \frac{\rho}{A^2}(W_1\bw_2 - \bw_1
W_2)]\d \xi\\
 & \quad + [\frac{1}{A}(\bw_2 \bp W_1 - W_1 \bp \bw_2 -
\bw_1 \bp W_2 + W_2 \bp \bw_1) + \frac{\bar{\rho}}{A^2} (W_1\bw_2
- \bw_1
W_2)] \d \bar{\xi} \Big\},\\
\d X_7 &= [\frac{1}{A} (W_1 \p \bw_2 - \bw_2 \p W_1 + W_2 \p \bw_1
- \bw_1 \p W_2) + \frac{\rho}{A^2} (W_1 \bw_2 + \bw_1 W_2)] \d \xi \\
&\quad +[\frac{1}{A}(\bw_2 \bp W_1 - W_1 \bp \bw_2 + \bw_1 \bp W_2
- W_2 \bp \bw_1) + \frac{\bar{\rho}}{A^2} (W_1 \bw_2 + \bw_1
W_2)]\d \bar{\xi},\\
\d X_8 &= [\frac{1}{A}(\p \bw_2 - \p W_2) + \frac{\rho}{A^2}
(\bw_2+W_2)]\d \xi + [\frac{1}{A} (\bp W_2 - \bp \bw_2) +
\frac{\bar{\rho}}{A^2} (\bw_2 + W_2)] \d \bar{\xi}.
\end{split}
\end{equation}
From the diagonal entries of the matrix $\d X^2$ we obtain
\begin{equation}
\label{5_9}
\begin{split}
\d X_3 &= 2 \Big\{[\frac{1}{A} (W_1 \p \bw_1 - \bw_1 \p W_1) +
\frac{\rho}{A^2} |W_1|^2]\d \xi \\
& \quad +[\frac{1}{A} (\bw_1 \bp W_1 - W_1 \bp \bw_1) +
\frac{\bar{\rho}}{A^2} |W_1|^2] \d \bar{\xi}\Big\},\\
\d X_4 &= 2 \Big\{[\frac{1}{A} (W_2 \p \bw_2 - \bw_2 \p W_2) +
\frac{\rho}{A^2} |W_2|^2]\d \xi \\
& \quad + [\frac{1}{A} (\bw_2 \bp W_2 - W_2 \bp \bw_2) +
\frac{\bar{\rho}}{A^2} |W_2|^2] \d \bar{\xi} \Big\}.
\end{split}
\end{equation}

Note that by virtue of the conservation law \Ref{1_16}, the
$1$-forms \Ref{5_8} and \Ref{5_9} are the exact differentials
of real-valued functions.  The functions 
$X_j (\xi, \bar{\xi})$ \, $j=1,\ldots,8$ constitute the
coordinates of the radius vector
\begin{equation}
\vec{X}(\xi,\bar{\xi}) = (X_1(\xi,\bar{\xi}),\ldots, X_8(\xi,
\bar{\xi}))
\end{equation}
of a two-dimensional surface in $\mathbb{R}^8$.  Thus, if the
complex-valued functions $W_i$, \,$i=1,2$ correspond to any solution of
the $\C P^2$ sigma model \Ref{1_13}, then we can use the
generalized Weierstrass formulae \Ref{5_8} and \Ref{5_9} to
construct a two-dimensional surface in $\mathbb{R}^8$ uniquely
defined by this solution.

Let us note that in the limiting case when 
\begin{equation}
\label{1_24} W_i \rightarrow \frac{W}{\sqrt{2}}, \quad \quad i=1,2,
\end{equation}
or when we put $W_1=0$ or $W_2=0$, the Weierstrass formulae \Ref{5_8} and \Ref{5_9}  reduce to formulae \Ref{4_3} describing the immersion in the the $\C P^1$ case.
These limits characterize some properties of solutions of both systems 
\Ref{1_22'}  and \Ref{1_23}. 

Now, let us discuss some classes of surfaces immersed in $\mathbb{R}^8$ which can be determined directly by applying the Weierstrass representation \Ref{5_8} and \Ref{5_9}.

\noindent\\
{\bf Example 1.} As well known, the simplest case of solutions of the $\C P^2$ model is obtained by taking $W_i$ to be analytic. 
In this case $\bar{\partial W_i}=0$ and so many expressions in \Ref{5_8} and \Ref{5_9}  vanish. In fact we get, (with $c.c.$ denoting the complex conjugate) 
\begin{equation}
\begin{split}
\d X_1 &=  \p\Big\{\frac{W_1+\bw_1}{A}\Big\}\d\xi\  +\  c.c.,\qquad
\d X_2 = -i \p\Big\{\frac{W_1-\bw_1}{A}\Big\}\d\xi\  +\  c.c.,\\
\d X_3 & = 2 \p\Big\{\frac{\vert W_1\vert\sp2}{A}\Big\}\d\xi\  +\  c.c.,\qquad
\d X_4 = 2 \p\Big\{\frac{\vert W_2\vert\sp2}{A}\Big\}\d\xi\  +\  c.c..\\ 
\d X_5 & = -i \p\Big\{\frac{W_2-\bw_2}{A}\Big\}\d\xi\  +\  c.c.,\qquad
\d X_6 = -i \p\Big\{\frac{\bw_1 W_2-\bw_2 W_1}{A}\Big\}\d\xi\  +\  c.c.,\\
\d X_7 &=  \p\Big\{\frac{\bw_1 W_2+\bw_2 W_1}{A}\Big\}\d\xi\  +\  c.c.,\qquad
\d X_8 =  \p\Big\{\frac{W_2+\bw_2}{A}\Big\}\d\xi\  +\  c.c. 
\end{split}
\end{equation}

These expressions can be easily integrated giving us, up to overall constants that can be 
added to any $X_i$:
\begin{equation}
\begin{split}
 X_1 &=  \Big\{\frac{W_1+\bw_1}{A}\Big\},\qquad 
 X_2 = -i \Big\{\frac{W_1-\bw_1}{A}\Big\},\qquad
 X_3  = 2 \Big\{\frac{\vert W_1\vert\sp2}{A}\Big\},\\
 X_4 &= 2 \Big\{\frac{\vert W_2\vert\sp2}{A}\Big\}\qquad
 X_5 = -i \Big\{\frac{W_2-\bw_2}{A}\Big\},\qquad  
 X_6 = -i \Big\{\frac{\bw_1 W_2-\bw_2 W_1}{A}\Big\},\\ 
 X_7 & =  \Big\{\frac{\bw_1 W_2+\bw_2 W_1}{A}\Big\},\qquad 
 X_8 =  \Big\{\frac{W_2+\bw_2}{A}\Big\}.  
 \end{split}
\end{equation}

Note that in general we have a surface in $\mathbb{R}^8$.
Using \Ref{4_2} it is very easy to calculate the curvature as we know that
\begin{equation}
 g_{\bar\xi,\xi}= A^{-2}\{|\dot{W_1}|^2 + |\dot{W_2}|^2 + |W_1\dot{W_2} - W_2\dot{W_1}|^2\},\quad  g_{\xi,\xi} = g_{\bar\xi,\bar\xi} = 0,
\end{equation}
where the dot denotes the differentiation with respect to $\xi$.
Then the Gaussian curvature is given by 
\begin{equation}
\label{holo1}
K\,=\,-{2\over g_{\bar\xi,\xi}}\,\p \bar{\p}\ln{g_{\bar\xi,\xi}}.
\end{equation}

Now, by setting $W_2=0$ the above model is reduced to the $\C P^1$ case, with  
\begin{equation} 
X_1=\frac{W_1+\bar{W_1}}{1+|W_1|^2}, \quad X_2=i\frac{\bar{W}_1- {W_1}}{1+|W_1|^2}, \quad 
X_3=2\frac{|W_1|^2}{1+|W_1|^2}
\end{equation}
and the remaining components $X_i=0$ (for $i= 4,\ldots,7$).

Note that in this case our surface is the surface of an appropriately located sphere. To see this note that
\begin{equation}
X_3\,=\,1\,+\,\frac{1-|W_1|^2}{1+|W_1|^2}.
\end{equation}
Then we have
\begin{equation}
X_1^2\,+\,X_2^2 \,+\,(X_3-1)^2\,=\,1
\end{equation}
and so we see that all the points lie on the surface of a sphere of unit radius, centred at $(0,0,1)$. Of course, the number of times this surface is covered depends on the degree of
$W_1$, i.e. the topological charge of the map.
This is, of course, consistent with \Ref{holo1}, which gives a constant.

In the $\C{P}^2$ the situation is more complicated but also more can be said about the 
surface; i.e., for example, all points lie on the hyperellipsoid surface
\begin{equation}
X_1^2+X_2^2+X_3^2+X_4^2+X_5^2+2X_6^2+2X_7^2+X_8^2=2.
\end{equation}

However, the Gaussian curvature is not necessarily constant. To see this we use \Ref{holo1}
and consider  very specific fields, namely:
\begin{equation}
W_1\,=\,a\xi,\qquad W_2\,=\,\xi^2.
\end{equation}
Then 
\begin{equation}
 g_{\bar\xi,\xi}\,=\,\frac{a^2+4\vert \xi\vert^2+a^2\vert \xi\vert^4}{(1+a^2\vert \xi\vert^2+\vert \xi\vert^4)^2}
\end{equation}
and so
\begin{equation}
\p\bar{\p} g_{\bar\xi,\xi}\,=\,4\frac{(a^2+a^2\vert \xi\vert^4+a^4\vert \xi\vert^2)}{(a^2+4\vert \xi\vert^2+a^2\vert \xi\vert^4)^2}
\,-\,2\frac{a^2+a^2\vert \xi\vert^4+4\vert \xi\vert^2}{(1+a^2\vert \xi\vert^2+\vert \xi\vert^4)^2}.
\end{equation}
Then 
\begin{equation}
K=-4\,+\,8a^2\frac{(1+a^2\vert \xi\vert^2+\vert \xi\vert^4)^3}{(a^2+4\vert \xi\vert^2+a^2\vert \xi\vert^4)^3}
\end{equation}
which shows that the curvature is constant only when $a=0$ or $a=\sqrt{2}$, but not in other cases. Hence, in general, the surfaces do not have a constant Gaussian curvature.

{\bf Example 2.} The simple mixed solution obtained by choosing $g_1 = 1$, $g_2 = \xi$, $g_3 = \xi^2$ in the formulae \Ref{6_1}--\Ref{5_3}  gives us the following  
\begin{equation}
W_1 = \frac{-\bar{\xi}(1+2|\xi|^2)}{\xi(2+|\xi|^2)} \,,\quad W_2 = \frac{1-|\xi|^4}{\xi(2+|\xi|^2)}.
\end{equation}
The Weierstrass representation \Ref{5_8}, \Ref{5_9} can be integrated and it leads to the following expression for the immersion of our surface in $\mathbb{R}^8$ in polar coordinates $(r,\,\varphi)$ by 
\begin{equation}
\begin{split}
X_1(r,\varphi) &= \frac{-12r^4\cos{2\varphi}}{(1+r^2+r^4)(1+4r^2+r^4)},\\
X_2(r,\varphi) &= -\frac{12r^4\sin{2\varphi}}{(1+r^2+r^4)(1+4r^2+r^4)},\\
X_3(r,\varphi) &= \frac{-4(4r^6+6r^4+9r^2+2)}{(1+r^2+r^4)(1+4r^2+r^4)},\\
X_4(r,\varphi) &= \frac{12r^2(1+r^4)}{(1+r^2+r^4)(1+4r^2+r^4)},\\
X_5(r,\varphi) &= -\frac{2(r^8+7r^6-r^2-1)\sin{\varphi}}{r(1+r^2+r^4)(1+4r^2+r^4)},\\
X_6(r,\varphi) &= -\frac{-4(r^8-2r^6-4r^2-1)\sin{\varphi}}{r(1+r^2+r^4)(1+4r^2+r^4)},\\
X_7(r,\varphi) &= \frac{-4(r^8-2r^6-4r^2-1)\cos{\varphi}}{r(1+r^2+r^4)(1+4r^2+r^4)},\\
X_8(r,\varphi) &= \frac{-2(r^8+7r^6-r^2-1)\cos{\varphi}}{r(1+r^2+r^4)(1+4r^2+r^4)}.
\end{split}
\end{equation}
The curvatures can be calculated, but the expressions are rather involved, so we omit them here.

\noindent\\
{\bf Example 3.} Another interesting class of mixed solutions of the $\C{P}^2$ model is given by
\begin{equation}
W_1 = \frac{\xi+\bar{\xi}}{1-|\xi|^2} \,,\quad W_2 = \frac{\bar{\xi}-\xi}{1-|\xi|^2}.
\label{5_20}
\end{equation}
In this case the system of Euler--Lagrange equations \Ref{1_23} simplifies considerably and therefore the Weierstrass formulae \Ref{5_8}, \Ref{5_9} can be easily integrated. In polar coordinates, setting $r=e^\vartheta$, we obtain 
\begin{equation}
\begin{split}
X_2(\vartheta, \varphi) &= e^{-\vartheta}\tanh\vartheta\sin{\varphi},\quad
X_8(\vartheta, \varphi) = e^{-\vartheta}\tanh\vartheta\cos{\varphi},\\
X_7(\vartheta, \varphi) &= e^{-\vartheta} \sech\vartheta, \quad
 X_1 = X_3 = X_4 =X_5 = X_6 = 0.
\end{split}
\label{5_21}
\end{equation}
This describes a surface of revolution which is contained in a subspace of dimension $3$ (see Fig. 1). The first and second fundamental forms are 
\begin{equation}
\begin{split}
I &= {1\over r^2(1+r^2)^2}\left[{1\over r^2}(r^4+6r^2+1)\d r^2 + (r^2-1)^2\d \varphi^2\right]\mbox{,}\\
II &= {4\over (1+r^2)^2(r^4+6r^2+1)^{1/2}}\left[(r^2+3)\d r^2 + r^2(r^2-1)\d \varphi^2\right]\mbox{.}
\end{split}
\label{5_22}
\end{equation}
The Gaussian and mean curvature are
\begin{equation}
\begin{split}
K &= {16r^8(r^2+3)\over (r^4+6r^2+1)^2(r^2-1)}\mbox{,}\\
H &= {r^4(r^4+4r^2-1)\over (r^4+6r^2+1)^{3/2}(r^2-1)}\mbox{.}
\end{split}
\label{5_23}
\end{equation}
Since the curvatures are not constant, the surface cannot be obtained from the $\C P^1$ model. 

\section{Concluding remarks and prospects for future developments}
\label{sec:Concl}
There are reasons to expect that the association
between the Weierstrass representation of surfaces immersed in
$\su{N+1}$ Lie algebras and the solutions of the Euclidean two-dimensional $CP^N$
sigma models, described in this paper, can be found also for more
general sigma models.  A good object of the investigation in this
direction are the complex Grassmannian sigma models which take
values on symmetric spaces $SU(m+n)/(S\left(U(m) \times
U(n)\right)$. These models share many important common properties
with the $CP^N$ models considered here.  They possess an infinite
number of conserved quantities, as well as infinite-dimensional
dynamical symmetries which generate the Kac-Moody algebra.  The
Grassmannian sigma model equations, just like those of the $\cpn$ models, have
a Hamiltonian structure and complete integrability with a well formulated linear
spectral problem.  Many classes of solutions of these
equations are known, see eg.\cite{WZ,GS2}; they can be expressed in terms of
holomorphic functions and functions obtained from them by 
a procedure which is a generalization of the
transformation which generates all solutions of the $\cpn$ models.

The complex Grassmannian sigma model in two Euclidean
dimensions are defined in terms of fields 
\begin{equation}
g = g(\xi,\,\bar\xi) \in \gr{SU}{N+1}
\end{equation}
where $\xi=\xi^1+i\xi^2$, taking values in the complex Grassmann manifold  
$\gr{SU}{N+1}/\mathbf{S}(\gr{U}{m}\times\gr{U}{n})$, where 
$N+1 = m+n$.  By decomposing $g$ into two blocks
$$
g=(X,Y), \, X=(z_1,\ldots, z_m), \, Y=(z_{m+1}, \ldots ,
z_{N+1}),
$$ 
where $z_i$ are $(N+1)-$component orthonormal column vectors,
\begin{equation}
z_i^{\dagger} \cdot z_k = \delta_{ik}
\end{equation}
we define the projector matrix $P$ ($P^{\dagger} = P$,  $P^2 = P$) as
\begin{equation}
P=X\,X^{\dagger} = \sum_{l=1}^m z_l \, z_l^{\dagger}.
\end{equation}
In general, it has higher rank than the corresponding matrix
for the $\cpn$ model. 
However, the equation of motion in terms of $P$ in this case has
the same form as \Ref{1_77} and is obtained by minimizing the
action of the Lagrangian
\begin{equation}
L=\mathrm{tr}\left[\left(D_{\mu} X\right)^{\dagger} \cdot
\left(D_{\mu} X \right)\right], \end{equation}
 where $D_{\mu} X$ is the covariant derivative for $X$,
\begin{equation}
D_{\mu} X = \p_{\mu} X -X (X^{\dagger}\cdot \p_{\mu} X).
\end{equation}
The above fact implies that our method can be successfully used
for constructing surfaces associated with the complex Grassmannian
sigma models.  The question of the diversity and complexity of
these surfaces, however, remains open and has to be answered in
further work.

In this paper we have shown how to generalize the old idea of
Enneper  \cite{Enne} and Weierstrass \cite{Weier}  in connection
with the $\cpn$ sigma models and their group properties. We have
found the structural equations of surfaces immersed in $\su{N+1}$ 
Lie algebras and expressed them in terms of any solution of the
$\cpn$ model. The most important advantage  of the presented method is that it is quite general. In constructing surfaces we proceeded directly from the given $\cpn$ model, without refering to any additional considerations. Another important advantage of our method is that, due to the conservation laws of the $\cpn$ model, the obtained expressions for surfaces are given at least in the form of quadratures.  

We have discussed in detail the geometrical aspects of the constructed surfaces. Namely, we have demonstrated through the use of moving frame that one can derive, via the $\cpn$ models, the first and the second  fundamental forms of a given surface as well as relations between them as expressed in the Gauss--Weingarten and the Gauss--Codazzi--Ricci equations. We have illustrated the proposed method of constructing surfaces in the case of low dimensional $\su{N+1}$ Lie algebras. 

A systematic application of the group theory makes its possible to obtain large number of particular solutions of the $\cpn$ equations and associated surfaces in $\NR^{N(N+2)}$. A question arises whether these solutions, corresponding to specific boundary conditions, are actually observable in nature. The answer depends to a large degree on their stability. Stable solutions should be observable and should also provide the starting point for perturbative calculations. These should in turn provide good approximative solutions relevant for situations in which the group-theoretical solutions no longer apply. In this context another question arises, namely what physical insight one gains from exact analytic expressions for surfaces. A particular answer is that they show up qualitative features that might be  difficult to detect numerically. We hope that our approach and our results may be useful in  applications to the study of surfaces which arise in physics, chemistry and biology, by providing explicit models in situations which have been well investigated experimentally but for which the theory is not yet well developed. Further exploration of relations between various properties of harmonic maps $S^2\to \cpn$ and properties of surfaces in planed for future work.

\noindent {\bf Acknowledgments} \\
Partial support for A.M.G. work was provided by a research grant
from NSERC of Canada and Fonds FQRNT du Gouvernement du Qu\'{e}bec.
A.S. wishes to acknowledge and thank the Centre de Recherches
Math\'{e}matiques (Universit{\'e} de Montr{\'e}al) and NSERC of Canada
for the financial support provided for his visit to Montreal.

\begin{footnotesize}

\end{footnotesize}

\label{lastpage}

\clearpage

\begin{figure}
\centerline{\epsfig{file=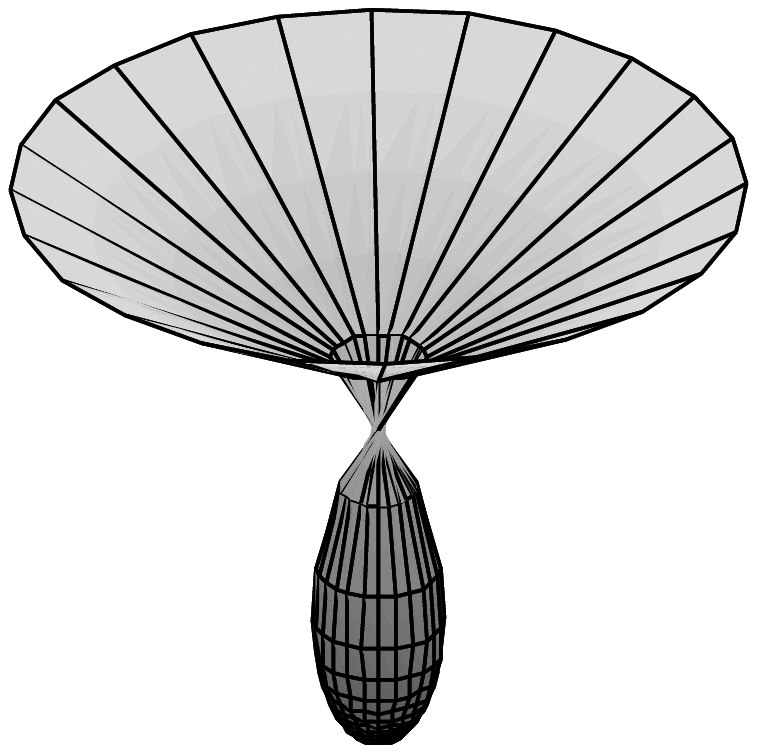,height=7.in} }
\caption{The surface associated with the solution \Ref{5_21}.}
\end{figure}

\clearpage

\end{document}